\documentclass[final]{siamltex}
\usepackage{amsmath}
\usepackage{amsfonts}
\usepackage{amssymb}
\usepackage[pdftex]{graphicx}
\usepackage{url}
\usepackage{color}
\usepackage{bm}

\newcommand{\om}{\omega}

\newcommand{\la}{\lambda}
\newcommand{\ep}{\varepsilon}

\newcommand{\oL}{\overline{L}}
\newcommand{\oT}{\overline{T}}

\newcommand{\cB}{\mathcal{B}}
\newcommand{\cC}{\mathcal{C}}
\newcommand{\cW}{\mathcal{W}}

\newcommand{\cU}{\mathcal{U}}
\newcommand{\cF}{\mathcal{F}}
\newcommand{\cR}{\mathcal{R}}
\newcommand{\cK}{\mathcal{K}}
\newcommand{\cD}{\mathcal{D}}

\newcommand{\cLe}{\mathcal{L}_\epsilon}
\newcommand{\cLo}{\overline{\mathcal{L}}}
\newcommand{\cT}{\mathcal{T}}
\newcommand{\cS}{\mathcal{S}}

\newcommand{\cJ}{\mathcal J}
\newcommand{\cI}{\mathcal I}
\newcommand{\cE}{\mathcal E}

\newcommand{\bu}{\mbox{u}}
\newcommand{\bp}{\mbox{p}}

\newcommand{\real}{\mathbb{R}}

\newcommand{\EE}{\mathbb{E}}
\newcommand{\PP}{\mathbb{P}}

\begin{document}
\title{Pulse propagation in time dependent randomly layered media}
\author{Liliana Borcea\footnotemark[1] \and Knut
  S{\O}lna\footnotemark[2]}
\renewcommand{\thefootnote}{\fnsymbol{footnote}}
\footnotetext[1]{Mathematics, University of Michigan, Ann Arbor, MI
  48109. (borcea@umich.edu)} \footnotetext[2]{Mathematics,
  University of California at Irvine, Irvine CA
  92697. (ksolna@math.uci.edu)} \maketitle \date{today}
\begin{abstract}
We study cumulative scattering effects on wave front propagation in
time dependent randomly layered media. It is well known that the wave
front has a deterministic characterization in time independent media,
aside from a small random shift in the travel time. That is, the pulse
shape is predictable, but faded and smeared as described
mathematically by a convolution kernel determined by the
autocorrelation of the random fluctuations of the wave speed.  The
main result of this paper is the extension of the pulse stabilization
results to time dependent randomly layered media. When the media
change slowly, on time scales that are longer than the pulse width and
the time it takes the waves to traverse a correlation length, the
pulse is not affected by the time fluctuations. {In rapidly changing
media, where these time scales are similar, both the pulse shape and
the random component of the arrival time are affected by the
statistics of the time fluctuations of the wave speed. We obtain an
integral equation for the wave front, that is more complicated than in
time independent media, and cannot be solved analytically, in general.
We also give examples of media where the equation simplifies, and the
wave front can be analyzed explicitly. We illustrate with these
examples how the time fluctuations feed energy into the
pulse. 
}
\end{abstract}

\section{{Introduction}}
\label{sect:intro} 
We study wave front propagation in time dependent, randomly layered
media. The initial condition corresponds to a plane wave normally
incident to the layers, so that the mathematical model for wave
propagation is the one dimensional wave equation. Extensions to three
dimensional wave fields generated by spatially localized sources are
not considered here, but may be done with plane wave decompositions
that lead to a family of one dimensional problems, as explained in
\cite[Chapter14]{newbible2007}. The medium fluctuations are around a
reference speed $c_o$ and are modeled with a random process $\nu(t,z)$
that is statistically stationary in time $t$ and homogeneous in range
$z$, that is, the propagation direction.  The problem is to describe
how a pulse impinging on the layered medium changes due to cumulative
scattering effects. The answer depends on the strength of the wave
speed fluctuations and the relation between the fundamental length and
time scales in the problem: The distance of propagation $L$, the
typical wavelength $\la_o$, the correlation length $\ell$ of the
layers, the localization length $L_{loc}$, and the correlation time
$T_\nu$ that quantifies the lifespan of a given spatial realization of
$\nu$. Localization means that the waves cannot penetrate too deep in
the random medium, and the incident energy is scattered back.  In
layered media a pulse is being transformed as it propagates a distance
that is on the scale of the localization length as described by the
pulse stabilization theory in
\cite{GeorgeP91Bible,burridge1988one,Clouet94,Solna00RayTheory} and in
this paper.
 
There are two high frequency regimes ($\la_o \ll L$) that have been
studied extensively in the context of wave propagation in time
independent randomly layered media (i.e., $\nu = \nu(z)$) because they
capture wave scattering effects in a canonical way
\cite{GeorgeP91Bible,burridge1988one,newbible2007}. They assume
distances of propagation comparable to the localization length ($L
\sim L_{loc}$) so that there is significant energy transfer between
the forward going and backscattered waves, and the regimes differ by
the strength of the fluctuations and the relation between the
wavelength and the correlation length \cite[Section
  5.1.4]{newbible2007}.  In the {\em weakly heterogeneous} regime the
wave speed fluctuations are small and $\ell \sim \la_o$, whereas in
the {\em strongly heterogeneous} regime the fluctuations are of order
one and $\ell \ll \la_o$. Interestingly, the wave front has a similar,
nearly deterministic characterization in both regimes, as shown in
\cite{GeorgeP91Bible,burridge1988one,Clouet94,Solna00RayTheory}. The
travel time is only slightly affected by the fluctuations, via a small
random shift defined by a standard Brownian motion, and the pulse
shape is a faded and smeared version of the initial one, described by
a deterministic convolution kernel determined by the autocorrelation
of $\nu(z)$.  This is called {\em pulse stabilization} and was first
noted by O'Doherty and Anstey in a geophysical context \cite{ODA71}.
A generalization of the pulse stabilization phenomenon persists in
media with long range correlations and in multi-scale media
\cite{GSfract,solna03}.

We consider time dependent randomly layered media. The scaling regime
is the weakly heterogeneous regime.  The waves interact strongly with
the medium in this regime, because the wavelength is on the scale of
the medium fluctuations, and the behavior of the pulse depends on the
full autocorrelation function of the fluctuations.  In the strongly
heterogeneous regime the waves do not see the layers in detail,
because $\la_o \gg \ell$, and the fluctuations take the ``effective''
form of white noise, independent of the detailed structure of the
random process $\nu$. The main result of the paper is the
characterization of the wave front in time dependent media. We show
that as long as the correlation time $T_\nu$ is long in comparison
with the time $T_\ell = \ell/c_o$ of travel over a correlation length,
the time changes in the medium have no effect on the pulse. The pulse
perceives the medium as if it were time independent.  In rapidly
changing media, where $T_\nu \sim T_\ell$, the wave front has a more
complicated behavior, since the medium is being transformed during the
passage time of the pulse.  We derive an integral equation for the
wave front. It has a deterministic kernel shape, however with a random
dilation, reflecting the random medium modulation during the pulse
passage time.  This equation is no longer solvable analytically, in
general. Nevertheless, there are examples of rapidly changing media
where the equation simplifies and the wave front has an explicit and
deterministic expression. We illustrate with such examples how the
temporal fluctuations feed energy into the wave front, thus enhancing
the pulse. That is to say, there is a trade-off between pulse fading
due to the spatial fluctuations and pulse enhancement due to the rapid
time changes in the medium.

Wave propagation in time dependent randomly layered media has also
been considered in \cite{frankenthal2004backscattering}, for
correlation times $T_\nu$ that are small with respect to the mean
scattering time\footnote{The mean scattering time is defined in
\cite[Equation (1.10)]{frankenthal2004backscattering}. It is an order
one time scale in our regime.}, but larger than
$T_\ell$. The formal stochastic iterative approach in
\cite{frankenthal2004backscattering} gives closed equations for the
statistical moments of the forward and backward going wave
fields. However, it relies on the assumption that the wave fields are
statistically independent from the medium fluctuations in certain
range intervals, which is difficult to justify rigorously. 

The paper is organized as follows: We begin in section
\ref{sect:formulation} with the mathematical formulation of the
problem and the asymptotic scaling regime. The statement of the pulse
stabilization result is in section \ref{sect:statement}, and its proof
is in section \ref{sect:proof}. We end with a summary in section
\ref{sect:summary}.

\section{Formulation of the problem}
\label{sect:formulation}
We consider the same one dimensional acoustic wave equation in time
changing media as in 
\cite{frankenthal2004backscattering}
\begin{eqnarray}
\rho \frac{\partial \bu (t,z)}{\partial t} + \frac{\partial
  \bp(t,z)}{\partial z} &=& F(t,z), \nonumber \\ \frac{1}{K(t,z)}
\frac{\partial \bp (t,z)}{\partial t} + \frac{\partial
  \bu(t,z)}{\partial z} &=& 0.
\label{eq:form.1}
\end{eqnarray}
It can be derived from the linearization of the mass and linear
momentum conservation equations in fluid dynamics, with $\bu$ the
displacement of the particles in the medium and $\bp$ the time
integral of the acoustic pressure.  The bulk modulus $K$ varies with
time $t$ and range $z$. The density $\rho$ is assumed constant for
simplicity, but the results can be extended to variable $\rho$.  We
model the variations of $K$, and thus of the wave speed $c =
\sqrt{K/\rho}$, with the random process $\nu(t,z)$
\begin{equation}
\frac{1}{c(t,z)} = \left\{ \begin{array}{l} \left[ 1 +
    \nu(t,z) \right]/c_o, \quad \mbox{ for} ~
     z> 0,
  \\ 1/c_o \quad \mbox{ otherwise}, \end{array} \right. 
\label{eq:form.2}
\end{equation}
and take for simplicity a constant reference wave speed $c_o$.  At any
time instant $ t \ge 0$ the fluctuations are confined to the half line
$z > 0.$

The medium is quiescent for $t < 0$
\begin{equation}
p(t,z) = 0, \qquad u(t,z) = 0, \quad t < 0,
\end{equation} 
and the waves are generated by a source located at 
$z = 0$ emitting a
pulse $\cF(t)$, 
\begin{equation}
F(t,z) = \zeta_o^{1/2} \cF(t) \delta(z/L),
\label{eq:form.3}
\end{equation}
for $L$ the range scale of propagation.  Here $ \zeta_o = \rho c_o $
is the reference acoustic impedance, and we normalize
(\ref{eq:form.3}) with it and $L$ to simplify the formulas below.  The
source generates a left (backward) going wave that never interacts
with the layered medium, and a right (forward) going wave which
penetrates the medium.  The problem is to describe the wave front,
that is the right going wave observed near the expected travel
time $\tau = z/c_o$ for $z > 0$, on a time scale comparable to the
width of the pulse $\cF(t)$.

\subsection{Random model of the fluctuations}\label{sec:medium} 
The random process $\nu(t,z)$ is a mathematical model of small scale
fluctuations of the wave speed. Small scale means that $c(t,z)$ 
varies rapidly in $z$ on the range scale $L$ of propagation of
the waves, and it changes over time intervals that are smaller than
the travel time scale $T = L/c_o$. The fluctuations are not known in
detail, as is typical in applications, and the random process
$\nu(t,z)$ models their uncertainty.
 
We take $\nu$ to be a mean zero stationary random field with standard
deviation $\sigma$
\begin{equation}
\EE \left\{ \nu(t,z) \right\} = 0, \qquad \sigma = \sqrt{\EE\left\{
  \nu^2(t,z) \right\}},
\label{eq:nu.2}
\end{equation}
and assume that it is bounded, so that the right hand
side in (\ref{eq:form.2}) remains positive. Moreover, $\nu$ is twice
continuously differentiable with   bounded first and
second order derivatives, and with rapid decay of 
correlations\footnote{See equation (\ref{eq:asMix}) for a more precise
  statement of the mixing assumption on the fluctuations.}. We scale
its range dependence by the correlation length $\ell$, and its time
dependence by the correlation time $T_\nu$.  These are related to the
covariance of $\nu$ as
\begin{align}
\hspace{-0.1in}\ell = \frac{1}{\sigma^2} \int_{-\infty}^\infty dz' \EE \left\{
\nu(t,z) \nu(t,z + z') \right\}, \quad 
\label{eq:nu.3} 
T_\nu =\frac{1}{\sigma^2} \int_{-\infty}^\infty dt' \EE \left\{ \nu(t,z)
\nu(t+t',z ) \right\}.
\end{align}

\subsection{Scaling}
\label{sect:formulation.1}
We model the pulse $\cF(t)$ using a real valued function $f(s)$
of dimensionless argument
\begin{equation}
\cF(t) = f \left({t}/{{T_{_\cF}}}\right),
\label{eq:pulse.2}
\end{equation}
and denote by $T_{\cF}$ the time scale that quantifies the pulse
duration.  We assume that $f$ is continuously differentiable and
compactly supported in the interval $[0,S]$, with Fourier transform
$\hat f(w)$ that peaks at $|w| = 2 \pi$ and is supported away from the origin.  
Since $f$ is real valued,
$\hat f(-w)$ is defined by the complex conjugate of $\hat f(w)$.  The
Fourier transform of the pulse
\begin{equation}
\widehat \cF(\om) = \int_\infty^\infty dt \,
f\left(\frac{t}{{T_{_\cF}}}\right) e^{i \om t} = T_{_\cF}
\widehat{f}\left(\om T_{_\cF}\right)
\end{equation}
has the central frequency $\om_o = 2 \pi/T_{_\cF}$, which gives the
central wavelength $\la_o = c_o T_{_\cF}$.

Let us rewrite the fluctuations in scaled form
\begin{equation}
\nu(t,z) = \sigma \mu \left(\frac{t}{T_\nu},\frac{z}{\ell} \right),
\label{eq:nu.5}
\end{equation}
using the function $\mu$ of dimensionless arguments. Its
autocorrelation is given by
\begin{equation}
\Phi(t',z') = \EE\left\{\mu(t''+t',z''+z') \mu(t'',z'') \right\},
\label{eq:nu.6}
\end{equation} 
for dimensionless $t',z' \in \real$. It is integrable and normalized
by
\begin{equation}
\Phi(0,0) = 1, \qquad \int_{-\infty}^\infty d t' \Phi(t',0) = 1, 
\qquad 
\int_{-\infty}^\infty d z' \Phi(0,z') = 1,
\label{eq:nu.7}
\end{equation}
to be consistent with definitions (\ref{eq:nu.2})-(\ref{eq:nu.3}).

The reference range scale is $\overline{L}$, the typical distance of
propagation through the medium, and the time is scaled by
$\overline{T} = \overline{L}/c_o$. We assume that we observe the waves for a
maximum time $T \sim \oT$, where ``$\sim$'' denotes equal up to a
dimensionless order one constant, and use causality to truncate
mathematically the fluctuations for $z > L$. This is because until
time $T$ the waves are not affected by the medium beyond the range
$L \approx c_o T$.

In the scaled coordinates $ t'= {t}/{\oT}$ and $z'= {z}/{\oL}$, the
pulse becomes
\begin{equation}
\cF(t) = f \left( \frac{t'}{T_{_\cF}/\oT}\right),
\label{eq:sc.1}
\end{equation}
and  the model of the wave speed is  
\begin{equation}
\frac{c_o}{c(t,z)} = \left\{ \begin{array}{l} \left[ 1 + \sigma \mu
    \left( \frac{t'}{T_\nu/\oT}, \frac{z'}{\ell/\oL} \right)\right]/c_o', \quad
  \mbox{ for} ~ ~z' > 0, \quad z' \le L'= {L}/{\oL} \sim 1,
  \\ {1}/{c_o'} \quad \mbox{otherwise}. \end{array} \right.
\label{eq:sc.2}
\end{equation}
In our case $c_o' = 1$ because the reference speed is constant, but in
general it would be a dimensionless function of $z'$. We keep $c_o'$
in the equations to distinguish between the range and time variables,
even though $c_o' = 1$.

\subsection{The asymptotic regime}
\label{sect:formulation.2}
We assume from now on that the time, range and speed are scaled and
drop the primes.  Our analysis is in an asymptotic regime modeled with
the small parameter $\ep$ defined as
\begin{equation}
\ep^2 = {T_\cF}/{\oT} = {\la_o}/{\oL} \ll 1.
\label{eq:as.3}
\end{equation}
It corresponds to a long distance of propagation compared to the
wavelength.  The fluctuations of the wave speed are weak, with
standard deviation\footnote{The standard deviation does not have to be
  exactly $\ep$, but of order $\ep$, that is $\sigma = \sigma_o \ep$,
  for some $\sigma_o = O(1)$. For simplicity we set $\sigma_o = 1$. A
  similar simplification is made in equation (\ref{eq:as.3}).}$ \sigma
= \ep,$ but they have a significant net scattering effect when the
correlation length $\ell$ is similar to $\la_o$, or the correlation
time $T_\nu$ is similar to the pulse width $T_{\cF}$. For example, it
is shown in \cite[Chapter 14]{newbible2007} that in time independent
media with $\ell= \la_o$, the localization length \cite[Section
  5.1.4]{newbible2007} $L_{loc}$ is similar to our reference order one
length scale $\oL$,
\begin{equation}
\frac{L_{loc}}{\oL} \sim \frac{1}{\sigma^2 \oL (2\pi/\lambda_o)^2 \ell} =
\frac{1}{4 \pi^2}.
\label{eq:Lloc}
\end{equation}
This means that the net scattering is so strong that the waves cannot
penetrate much deeper than the range $\oL$ in the medium, and are
reflected back toward the source.

We model the range and temporal scales of the fluctuations by 
\begin{equation}
\ep^\beta = {\ell}/{\oL}, \qquad \ep^\alpha = {T_{\nu}}/{\oT}, 
\label{eq:as.1}
\end{equation}
using the dimensionless parameters $\alpha$ and $\beta$ in the
interval $[0,2]$. When $\alpha <2$ we say that the medium is slowly
changing, because it varies little during the duration of the pulse.
In rapidly changing media $\alpha = 2$. Similarly, $\beta< 2$ models
a medium that varies slowly in $z$.

In our scaling the pulse becomes $f\left(t/\ep^2\right)$ and the model
of the wave speed in the random medium is
\begin{equation}
c_\ep(t,z) = {c_o} \Big[ 1 + \ep
\mu\Big(\frac{t}{\ep^\alpha},\frac{z}{\ep^\beta}\Big)
\Big]^{-1},
\label{eq:as.5}
\end{equation}
where we use the index $\ep$ to remind us of the dependence on the
asymptotic parameter $\ep \ll 1.$ The acoustic impedance is given by
the product of the wave speed and the constant medium density, so it
takes the form
\begin{equation}
\zeta_\ep(t,z) = {\zeta_o}{ \Big[ 1 + \ep
\mu\Big(\frac{t}{\ep^\alpha},\frac{z}{\ep^\beta}\Big)
\Big]^{-1}}.
\label{eq:as.6}
\end{equation}
 
\section{Statement of results}
\label{sect:statement}
The goal of the paper is to characterize the wave front in the
asymptotic limit $\ep \to 0$.  To define it we introduce the right and
left going waves $A_\ep(t,z)$ and $B_\ep(t,z)$, moving with the local
speed $c_\ep(t,z)$, as in \cite[Section 8.1]{newbible2007}
\begin{eqnarray}
A_\ep(t,z) &=& \zeta_\ep^{-1/2}(t,z) \bp_\ep(t,z) + \zeta_\ep^{1/2}(t,z)
\bu_\ep(t,z), \nonumber \\ B_\ep(t,z) &=& -\zeta_\ep^{-1/2}(t,z)
\bp_\ep(t,z) + \zeta_\ep^{1/2}(t,z) \bu_\ep(t,z).
\label{eq:dec.1}
\end{eqnarray}
Here $u({\oT} t, {\oL} z) = u_\ep(t,z)$, with index $\ep$ indicating the
dependence on the asymptotic parameter $\ep$. The same notation
applies to the pressure field $\bp_\ep(t,z)$.  The right going wave
equals the impinging pulse at $z=0$,
\begin{equation}
A_\ep(t,0) = f \left({t}/{\ep^2}\right).
\label{eq:dec.2}
\end{equation}
At the order one scaled final range $L$ we impose the condition
\begin{equation}
B_\ep(t,L) = 0.
\label{eq:dec.3}  
\end{equation} 
It follows from the truncation of the random medium in
(\ref{eq:sc.1}), justified by the causality of the wave equation, as
explained in section \ref{sect:formulation.1}. It says that up to time 
$T$ we cannot observe waves backscattered beyond the range $L$.

\subsection{The wave front and the random travel time}
\label{sect:WFRT} 
We define the wave front at range $z$ as the right going wave observed
in an $\ep^2$  time window around  the reference travel time $\tau = z/c_o$,
\begin{equation}
a_\ep(\tau,s) = A_\ep(t_\ep(c_o \tau,s),c_o\tau).
\label{eq:wf}
\end{equation} 
Here 
\begin{align}
t_\ep(c_o \tau,s) &= \ep^2 s + \int_0^{c_o \tau} \hspace{-0.1in}
\frac{dz'}{c_\ep(t_\ep(z',s),z')} = \tau + \ep^2 s + \ep^2
\cW_\ep(\tau,s),
\label{eq:tt}
\end{align} 
is the travel time along the characteristics of the right going wave,
as explained in section \ref{sect:proof.char}. The characteristics are
parameterized by the range $z = c_o \tau$, and start at times $\ep^2 s$
at $z = 0$, with $s$ parameterizing the pulse. Note that $t_\ep(c_o
\tau,s)$ fluctuates randomly around the value $\tau + \ep^2 s$, on the
scale of the duration of the pulse, as modeled by the process
\begin{equation}
\cW_\ep(\tau,s) = \frac{1}{\ep} \int_0^\tau d u \, \mu \left(
\frac{u}{\ep^\alpha} + \ep^{2-\alpha}(s +
\cW_\ep(u,s)),\frac{c_o u}{\ep^\beta}\right).
\label{eq:Rshift}
\end{equation} 

The wave front (\ref{eq:wf}) satisfies the following equation, as
proved in section \ref{sect:proof.inteq},
\begin{align} 
a_\ep(\tau,s)  =& f(s) - \frac{c_o^2}{8} \int_0^\tau d \tau'
\mu_{\alpha,\beta}^+ \Big[\frac{\tau' + \ep^2\left(s +
    \cW_\ep(\tau',s)\right)}{\ep^\alpha}, \frac{c_o \tau'}{\ep^\beta}
  \Big] \int_{Y^\ep_{\alpha,\beta}(\tau',s)}^s dy\, \times 
\nonumber  \\
&\mu_{\alpha,\beta}^- \Big[\frac{ \tau' + \ep^2 \left(\frac{y+s}{2}+
    \cW_\ep(\tau',s)\right)}{\ep^\alpha}, \frac{c_o \tau' + \ep^2 c_o
    (s-y)/2}{\ep^\beta}\Big] \times \nonumber \\
    &a_\ep
\left(\tau'+ \ep^2(s-y)/2,\cS_\ep(y;\tau',s)\right) + O(\ep), 
\label{eq:inteq} 
\end{align}
where the integration bound $Y_{\alpha,\beta}^\ep(\tau',s) < s$ is given explicitly
in section \ref{sect:proof.inteq}.    The random
processes $\mu^\pm_{\alpha,\beta}$ are
\begin{equation}
\mu_{\alpha,\beta}^\pm (t,z) = \ep^{2-\beta} \mu_z(t,z) \pm
{\ep^{2-\alpha}}c_o^{-1} \mu_t(t,z),
\label{eq:nupm}
\end{equation}
with $\mu_z$ and $\mu_t$ the partial derivatives in $z$ and $t$ of
$\mu$. The mapping $\cS_\ep(y;\tau,s)$ is defined by the equation
\begin{align}
\cS_\ep(y;\tau,s) + \cW_\ep(\tau, \cS_\ep(y;\tau,s)) =
y + \cW_\ep(\tau,s) +  \ep
\cR_\ep(y,\tau,s),\label{eq:mapSe}
\end{align}
where $\cR_\ep$ is bounded\footnote{The process
  $\cR_\ep$ is defined explicitly later, in section
  \ref{sect:proof.inteq}.}.

The complication in equation (\ref{eq:inteq}) stems from the possible
dependence of the random fluctuations (\ref{eq:Rshift}) of the travel
time on $s$, the parameter of the pulse. This is why we need the
mapping $\cS_\ep$, and the integration bound $Y_{\alpha,\beta}^\ep$.
To obtain an explicit characterization of the wave front, we consider
two scaling regimes:

\textbf{Regime 1} defined by $\alpha \le 1$ and $\beta =2$, where the 
medium is changing slowly.

\textbf{Regime 2} defined by $\alpha  = 2$ and $ \beta \le 1$, where the 
medium is changing rapidly.
\\
The random travel time fluctuations are independent of $s$ in these regimes
\footnote{The results can be extended to $\alpha < 2$ for regime 1 and
  $\beta < 2$ for regime 2, but the analysis is much more involved and
  will be presented in a different publication.}  as  stated in the next
lemma, proved in Appendix \ref{sect:Pf}.  The implication is that we
can set $ Y_{\alpha,\beta}^\ep = 0$ and $\cS_\ep(y;\tau,s)
\stackrel{\ep \to 0}{\longrightarrow} y$, where the convergence is in
probability. Equation (\ref{eq:inteq}) simplifies, and we can analyze
explicitly its solution to obtain the pulse stabilization results
described in the next section.

The analysis is based on the following mixing conditions satisfied by the process $\mu$. 
They can be derived from the  $\phi-$mixing  assumption on $\mu$, with mixing rate $(\Delta z)^{-d}$, as shown in appendix \ref{ap:MIXING}. 
There exists a deterministic constant $C$ and a real number $d > 2$ such that 
\begin{equation}
\EE\left[ \mu(t,z+\Delta z) \Big| \mathcal{F}_z\right] \le C (\Delta
z)^{-d}, 
\label{eq:asMix}
\end{equation}
for all $z, t$ and $\Delta z > 0$, where $\cF_z$ is the filtration associated to $\mu$. Moreover, suppose that $z_1$ and $z_2$ satisfy  
$|z_1-z| \ge \Delta z$ and $|z_2-z| \ge \Delta z$. Then, 
for all $z, t $ and $t'$ we have  
\begin{equation}
\left| \EE \left[ \mu( t,z_1) \mu(t',z_2 )
  \Big|\cF_{z}\right] - \EE\left[\mu (t,z_1) \mu(t',z_2 )
  \right] \right| \le C (\Delta z)^{-d}, 
\label{eq:asMix2}
\end{equation}
Similar conditions are assumed for the derivatives  $\mu_t$ and $\mu_z$ of the process, as well.   

\begin{lemma}
\label{thm.01}
Under the mixing assumption (\ref{eq:asMix}) for $\mu$, and similar for $\mu_t$, we have the following results:  In regime 1 \[
\cW_\ep(\tau,s) = W_\ep(\tau) + q_\ep(c_o \tau,s),
\]
where 
\begin{equation}
W_\ep(\tau) = \frac{1}{\ep} \int_0^\tau du \,
\mu\left(\frac{u}{\ep^\alpha},\frac{c_0 u}{\ep^2}\right),
\label{eq:defchi}
\end{equation} 
and  the residual $q_\ep$ satisfies
\begin{equation}
\PP\left(\sup_{s\in [0,S], \tau \in [0, {T}]} |q_\ep(c_o \tau,s)| \ge
c \ep^{(2 - \alpha - \theta)/2} \right) \le C \ep^\theta,
 \label{lem2.eq2}
\end{equation}
for $\ep$ independent constants $c$ and $C$, and any $\theta \in (0,
2-\alpha)$. Moreover,
\begin{equation} 
\PP\left(\sup_{s\in [0,S], \tau \in [0, {T}]} |\partial_s
\cW_\ep(\tau,s)| \ge c \ep^{(2 - \alpha - \theta)/2} \right) \le C
\ep^{\theta}.
 \label{lem2.eq1}
\end{equation}
 The result is similar in regime $2$, under the corresponding mixing assumptions on $\mu_z$,
 with $\alpha$ replaced by
 $\beta$ and  
 \begin{equation}
 \label{eq:reg2}
 W_\ep(\tau) = \frac{1}{\ep} \int_0^\tau du \,
\mu\left(\frac{u}{\ep^2},\frac{c_0 u}{\ep^\beta}\right) - \int_0^\tau du \,
\mu^2\left(\frac{u}{\ep^2},\frac{c_0 u}{\ep^\beta}\right).
\end{equation}
\end{lemma}

Thus, in regimes $1$ and $2$ the random travel time fluctuations are
described by the random process  $W_\ep(\tau)$. 
Under the mixing assumptions (\ref{eq:asMix}) and (\ref{eq:asMix2}), we   find \cite{heyde}  that   $W_\ep(\tau)$ in equation  (\ref{eq:defchi})  converges in distribution to 
the Gaussian process 
\begin{equation}
W(\tau) = h \tau + \cD \cB(\tau),
\label{eq:LIMW}
\end{equation}
where $\cB(\tau)$ is standard Brownian motion. The diffusion coefficient $\cD$ can be calculated  from the 
variance 
\begin{eqnarray}
   \lim_{\ep \to 0} 
   \EE \left[
   \frac{1}{\ep} \int_0^\tau du \,
\mu\left(\frac{u}{\ep^\alpha},\frac{c_0 u}{\ep^2}\right) \right]^2
=    \tau  \int_{-\infty}^\infty dt \, \Phi(0,c_ot)   =: \tau \cD^2
\label{eq:Mdiff}
\end{eqnarray}
in regime $1$, and 
\begin{eqnarray}
   \lim_{\ep \to 0} 
   \EE \left[
   \frac{1}{\ep} \int_0^\tau du \,
\mu\left(\frac{u}{\ep^2},\frac{c_0 u}{\ep^\beta}\right) \right]^2
=    \tau  \int_{-\infty}^\infty dt \, \Phi(t,0)   =: \tau \cD^2
\label{eq:Mdiffreg2}
\end{eqnarray}
in regime $2$.  The drift $h$ is obtained from the limit of the expectation of $W_\ep(\tau)$.  
It is $h = 0$ in regime $1$, and in regime $2$ we get $h = -1$, because  
\begin{equation}
\lim_{\ep \to 0} \EE \left[ W_\ep(\tau)\right] =- \lim_{\ep \to 0} \EE \left[ \int_0^\tau 
du \mu^2 \left(\frac{u}{\ep^2},\frac{c_0 u}{\ep^\beta}\right)\right] = - \tau .
\end{equation}

\subsection{Pulse stabilization}
\label{sect:PS}
The integral equation (\ref{eq:inteq}) of the wave front
$a_\ep(\tau,s)$ simplifies in regimes 1-2, with $Y_{\alpha,\beta}^\ep$
replaced by $0$ and $\cS_\ep(y;\tau,s)$ replaced by $y$. Moreover,
$a_\ep(\tau,s)$ converges to a deterministic function
$\bar{a}(\tau,s)$, as stated next.
\begin{theorem}
\label{thm.1} 
Assume regime $1$ or $2$. The wave front $a_\ep(\tau,s)$ converges in
probability as $\ep \to 0$, in the space of continuous functions of
$(\tau,s) \in [0,T]\times [0,S]$, to the deterministic solution
$\overline{a}(\tau,s)$ of the initial value problem
\begin{eqnarray}
\frac{\partial}{\partial \tau} \overline{a}(\tau,s) &=&
-\frac{c_o^2}{8} \int_{0}^s dy \, \overline{a}(\tau,y)
\Psi\left(\frac{s-y}{2} \right) , \quad \tau > 0, \nonumber
\\ \overline{a}(0,s) &=& f(s),
\label{eq:limit.ds}
\end{eqnarray}
with  integral kernel
\begin{equation}
\Psi(s) = \left\{\begin{array}{ll} 
-\partial_z^2 \Phi\left(0,c_o s\right) \quad &\mbox{in regime 1}, \\ \\
c_o^{-2}\partial_t^2 \Phi\left(s,0\right) & \mbox{in regime 2}.
\end{array} \right.
\label{eq:Psi1}
\end{equation}
The limit is given explicitly by 
\begin{equation}
\overline{a}(\tau,s) = \int_{-\infty}^\infty \frac{d \om}{2 \pi} \,
\hat f(\om) \exp\left[-i \om s - \frac{c_o^2 \tau }{4} \int_0^\infty
  du \, \Psi(u) e^{2 i \om u} \right].
\label{eq:PST}
\end{equation}
\end{theorem}

The convergence of $a_\ep$ to $\overline{a}$, the solution of
(\ref{eq:limit.ds}), is derived  in section \ref{sect:slow}. 
We assume stationary smooth medium fluctuations as explained in 
section \ref{sec:medium}, moreover, mixing conditions in (\ref{eq:asMix}) and (\ref{eq:asMix2})
in Regime 1 and correspondingly in Regime 2.

To obtain
(\ref{eq:PST}) note that the right hand side in (\ref{eq:limit.ds}) is
a convolution in $y$. Indeed, let us change variables in
(\ref{eq:limit.ds}) as $y = s-u$, for $u \in [0,s]$, and extend the
integral over $u$ to infinity
\begin{equation}
\frac{\partial}{\partial \tau} \overline{a}(\tau,s) = -\frac{c_o^2}{8}
\int_0^{\infty} du \, \overline{a}(\tau,s-u) \Psi\left(\frac{u}{2}
\right) ,
\label{eq:convol}
\end{equation}
because $\overline{a}(\tau,s-u) = 0$ for $s-u < 0$, by the compact
support of the pulse in $\mathbb{R}^+$. The Fourier transform
$\widehat{\overline{a}}(\tau,\om)$ with respect to $s$ satisfies
\[
\frac{\partial \widehat{\overline{a}}(\tau,\om)}{\partial \tau} = -
\frac{c_o^2}{4} \, \widehat{\overline{a}}(\tau,\om) \int_0^\infty du
\, \Psi(u) e^{2 i \om u}, \quad \tau > 0,
\]
with initial condition $\widehat{\overline{a}}(0,\om) = \widehat
f(\om).$ This can be solved explicitly and (\ref{eq:PST}) follows by
inverting the Fourier transform.

\subsection{Interpretation of the results}
\label{sect:discussion}
Theorem \ref{thm.1} states that when we observe the right going wave
at range $z$ around the random travel time $t_\ep (z,s)$, in a time
interval equal to the duration of the pulse, the result is
deterministic.  This phenomenon is called \emph{pulse stabilization.}
Note however that the pulse shape is not the same as the emitted
$f(s)$, like it would be in the homogeneous medium. The pulse is
deformed, as modeled by the convolution kernel
\begin{equation}
\mathcal{K}_f(z,s) = \int_{-\infty}^\infty \frac{d \om}{ 2 \pi} \,
e^{- i \om s - \om^2 \gamma(\om) z}, \qquad z = c_o \tau,
\label{eq:CONVKer}
\end{equation}
where we introduced the notation
\begin{equation}
\gamma(\om) = \frac{c_o}{4 \om^2} \int_0^\infty d u \, \Psi(u) e^{2 i
  \om u}.
\label{eq:ConvKer1}
\end{equation}

When we observe the wave at range $z$ in a time window centered at the
reference time $\tau = z/c_o$, the result is no longer deterministic,
due to the random component $\cW_\ep$ of the travel time. We state it
in the corollary below, which follows from Theorem \ref{thm.1}  
and  Slutsky's theorem: 
\begin{corollary}
\label{cor.1s}
As $\ep \to 0$,
\begin{equation}
A_\ep(\tau + \ep^2 s, c_o \tau) \to \int_{-\infty}^\infty \frac{d \om}{2
  \pi} \widehat f(\om) \, e^{ -i \om \left[ s -h \tau - \cD
    \cB(\tau) \right] - c_o \tau \om^2 \gamma(\om) },
\label{eq:corr_s}
\end{equation}
where the convergence is in distribution, $\cB(\tau)$ is standard
Brownian motion and $h$ and $\cD$ are as in (\ref{eq:LIMW}).
\end{corollary}

\subsubsection{Slowly time changing media}
 In regime 1 we have
\begin{align}
\rm{real}[\gamma(\om)] &= \frac{c_o}{4 \om^2} \int_0^\infty du \Psi(u)
\cos(2 \om u) = -\frac{1}{4 \om^2} \int_0^\infty dz\,
\partial_z^2 \Phi(0,z) \cos\left(\frac{2 \om}{c_o} z\right) \nonumber
\\ &= \frac{1}{2 c_o^2} \int_{-\infty}^\infty dz \, \Phi(0,z)
\cos\left(\frac{2 \om}{c_o} z\right) = \gamma^{(c)}(\om)
\ge 0,
\label{eq:Slow1}
\end{align}
where the second equality is from definition (\ref{eq:Psi1}), the
third equality is obtained with integration by parts, and the fourth
equality is a definition. That $\gamma^{(c)}$ is non-negative follows
by Bochner's theorem. The imaginary part of $\gamma(\om)$ is obtained
similarly
\begin{align}
\rm{imag}[\gamma(\om)] &= \frac{c_o}{4 \om^2} \int_0^\infty du \Psi(u)
\sin(2 \om u) = -\frac{1}{4 \om^2} \int_0^\infty dz\, \partial_z^2
\Phi(0,z) \sin\left(\frac{2 \om}{c_o} z\right) \nonumber \\ &=
\frac{1}{c_o^2} \int_{0}^\infty dz \, \Phi(0,z) \sin\left(\frac{2
  \om}{c_o} z\right) - \frac{1}{2 \om c_o} =\gamma^{(s)}(\om) -
\frac{1}{2 \om c_o}.
\label{eq:Slow2}
\end{align}

We recovered the pulse stabilization result \cite[Proposition
  8.1]{newbible2007} in time independent media. The random changes in
the medium, on the time scale $T_\nu = \ep^{\alpha}T$ that is much
larger than the pulse width $T_{_\cF}$ and the time $T_{\ell}$ of
traversal of a correlation length, play no role in the leading behavior of the
wave front. The result is as if the medium were time independent. The
kernel 
\begin{equation}
\mathcal{K}(z,s) = \int_{-\infty}^\infty \frac{d \om}{ 2 \pi} \, e^{-
  i \om s - z \om^2 [\gamma^{(c)}(\om) + i \gamma^{(s)}(\om)]} =
\mathcal{K}_f \left(z,s -\frac{z}{2c_o}\right),
\label{eq:CONVKer2}
\end{equation}
 is analyzed in detail in
\cite{newbible2007,GSfract}. It satisfies the initial value problem
\begin{equation}
\label{eq:pdeeff}
\frac{\partial \cK(z,s)}{\partial z} = {\cal L} \cK(z,s) \quad {\rm for} ~ z >
0, \quad {\rm and} ~ ~  \cK(z=0,s) = \delta(s),
\end{equation}
with pseudo-differential operator ${\mathcal L}$ describing the pulse
deformation caused by scattering in the medium. We can decompose it in
two parts, ${\mathcal L} = {\mathcal L}^{(c)} + {\mathcal L}^{(s)}$,
where
\begin{eqnarray}
\label{defl}
\int_{-\infty}^\infty {\mathcal L}^{(c)} \cK(z,s) e^{i \omega s} ds
&=& - \om^2 \gamma^{(c)}(\omega) \int_{-\infty}^\infty \cK(z, s) e^{i
  \omega s} ds \, , \\ \int_{-\infty}^\infty {\mathcal L}^{(s)} \cK(z,
s) e^{i \omega s} ds &=& - i \om^2 {\gamma^{(s)}} (\omega)
\int_{-\infty}^\infty \cK(z, s) e^{i \omega s} ds \, . \nonumber
\end{eqnarray}
The first part ${\mathcal L}^{(c)}$ acts as an effective diffusion
operator. It models the frequency-dependent attenuation over range
$z$, at the rate $\om^2 \gamma^{(c)}(\om)$.  The second part
${\mathcal L}^{(s)}$ is an effective dispersion operator that
preserves energy.  The dispersive effect increases with $z$ at rate
${\om^2 \gamma^{(s)}(\om)}$.  We can rewrite (\ref{defl}) in the time
domain using (\ref{eq:Slow1})-(\ref{eq:Slow2})
\begin{eqnarray}\nonumber
{\mathcal L} \cK(z,s) &=& \left[ \frac{1}{2 c_o} \Phi \left(0,
  \frac{c_o s}{2}\right) {\bf 1}_{[0,\infty)}(s) \right] * \left[
    \frac{\partial^2 \cK}{\partial s^2} (z,s)\right] \\ &=& \frac{1}{2
    c_o} \int_0^\infty \Phi \left(0, \frac{c_o \xi}{2}\right)
  \frac{\partial^2 \cK}{\partial s^2} (z, s-\xi) d\xi \, .
\end{eqnarray}
The indicator function ${\bf 1}_{[0,\infty)}$ entails that if $\cK$ is
  vanishing for $s< s_0$, then ${\mathcal L} \cK$ is also vanishing
  for $s < s_0$. That is to say, the operator ${\mathcal L}$ preserves
  causality. The operator also satisfies the Kramers-Kronig relations
  \cite{kronig,kramers}, as shown in \cite{GSfract}.

Corollary \ref{cor.1s} gives that $ A_\ep\left(\tau + \ep^2 s , c
\tau\right) \stackrel{\ep \to 0}{\longrightarrow} f_\tau \left(s -
\frac{\tau}{2} - \cD \cB(\tau) \right)$ in distribution, where $
f_\tau(s)= f(s) \star \cK\Big({c \tau }, s \Big)$ is the deterministic
pulse shape given by the convolution of $f(s)$ with the kernel
$\cK(z,s)$. The diffusion coefficient $\cD$ in front of the Brownian
motion is  
$ \cD = \sqrt{2 c_o {\gamma^{(c)}(0)}}.$ Note that the result is the
same as in \cite[Proposition 8.1]{newbible2007}, except for the
deterministic delay of $\tau/2$. This is because we have a different
model of the random wave speed. In \cite[Proposition
  8.1]{newbible2007} the denominator in (\ref{eq:as.5}) has a square
root, and this results in a deterministic correction of the travel
time that cancels the delay.

\subsubsection{Rapidly time changing media}
In regime 2 we have 
\begin{align}
\rm{real}[\gamma(\om)] &= \frac{c_o}{4 \om^2} \int_0^\infty du \Psi(u)
\cos(2 \om u) = \frac{1}{4 c_o\om^2} \int_0^\infty \, \partial_t^2
\Phi(t,0) \cos\left(2 \om t\right) \nonumber \\ &= -\frac{1}{2 c_o}
\int_{-\infty}^\infty dt \, \Phi(t,0) \cos\left(2 \om t\right) =
\gamma^{(c)}(\om) \le 0,
\label{eq:Fast1}
\end{align}
where the second equality is from definition (\ref{eq:Psi1}), the
third equality is obtained with integration by parts, and the fourth
equality is a definition. Note that the real part of $\gamma$ is now
negative. The imaginary part of $\gamma(\om)$ is obtained similarly
\begin{align}
\rm{imag}[\gamma(\om)] &= \frac{c_o}{4 \om^2} \int_0^\infty du \Psi(u)
\sin(2 \om u) = \frac{1}{4 c_o\om^2} \int_0^\infty dt\,
\partial_t^2 \Phi(t,0) \sin\left({2 \om} t\right) \nonumber
\\ &= -\frac{1}{c_o} \int_{0}^\infty dt \, \Phi(t,0)
\sin\left(2 \om t\right) + \frac{1}{2 \om c_o}
=\gamma^{(s)}(\om) + \frac{1}{2 \om c_o}.
\label{eq:Fast2}
\end{align}
As we did before, we define the pulse shaping kernel
$\mathcal{K}(z,s)$ by
\begin{equation}
\mathcal{K}(z,s) = \int_{-\infty}^\infty \frac{d \om}{ 2 \pi} \, e^{-
  i \om s - z \om^2 [\gamma^{(c)}(\om) + i \gamma^{(s)}(\om)]} =
\mathcal{K}_f \left(z,s{+}\frac{z}{2c_o}\right).
\label{eq:Fast3}
\end{equation}
This is to be understood weakly, as the finite bandwidth makes the
convolution $f_\tau(s) = f(s) \star_s \mathcal{K}(z,s)$ well defined.
We obtain  that
\begin{equation}
A_\ep\left(\tau + \ep^2 s , c \tau\right) \stackrel{\ep \to
  0}{\longrightarrow} f_\tau \left(s + \frac{3\tau}{2} - \cD \cB(\tau)
\right),
\label{eq:Fast4}
\end{equation} 
in distribution. Unlike the result in slowly changing media where the
pulse fades as it travels through the medium, in rapidly changing
media the pulse is enhanced. This is because $\gamma^{(c)}(\om)$ is
negative and energy is not conserved in time dependent media. Thus,
the random fluctuations of the wave speed in time and range have
opposite effects.  We illustrate this with a simple example in section
\ref{sect:EX_R}. Before that, we draw an analogy between the pulse
stabilization result in rapidly changing media and the random harmonic
oscillator problem studied in \cite[Section 7.5]{newbible2007} and \cite{KeP}.

\subsubsection{Analogy to the random harmonic oscillator problem}

Let us suppose for a moment that there were no spatial fluctuations,
just the temporal ones. We denote these fluctuations by $\widetilde
\mu$ to distinguish them from the time and range fluctuations in
(\ref{eq:as.5}).  We also change variables $ z = c_o {\cT}, t =
Z/c_o,$ and obtain a new wave equation for the displacement $\bu$,
\[
\frac{\partial ^2 \bu(\cT,Z)}{\partial Z^2} - \frac{\Big[ 1 + \ep
    \tilde \mu \Big( \frac{Z/c_o}{\ep^2} \Big) \Big]}{c_o^2}
\frac{\partial^2 \bu(\cT,Z)}{\partial \cT^2} = 0, \quad Z > 0,
\]
where \[ 1 + \ep \tilde \mu \left( \frac{Z/c_o}{\ep^2} \right) = \Big[
  1 + \ep \mu \Big( \frac{Z/c_o}{\ep^2}, 0 \Big) \Big]^{-2}.\]
Before we had a two point boundary value problem for the right and
left going waves. The right going wave amplitude was specified in
terms of the impinging pulse $f(t/\ep^2)$ by (\ref{eq:dec.2}) at $z =
0$. The left going wave amplitude was set to zero at range
$L$, using causality and the finite time $T$ of
observation. Now, after the change of variables, we have initial
conditions at $Z = 0$, which corresponds to $t = 0$,
\begin{align*}
\bu(\cT,Z =0  ) = f(\cT/\ep^2 )/(2\sqrt{\zeta_o}), \qquad 
\frac{\partial \bu(\cT,Z = 0 )}{\partial Z} = f'(\cT/\ep^2 )/( 2
\ep^{2} c_o \sqrt{\zeta_o}) .
\end{align*}
The first initial condition follows from the wave decomposition
(\ref{eq:dec.1}) evaluated at $Z = 0$ (i.e., $t=0$), where the wave is
right going with amplitude given by the pulse $f(\cT/\ep^2)$. The
second condition follows from the same decomposition and the relation
$\partial_{_Z} u(0,\cT) = -1/(c_o \zeta_0) \partial_{_\cT} p(0,\cT)$
obtained from the acoustic system of equations.

Fourier transforming in $\cT$, we obtain that $ \widehat
\bu_\ep(\om,Z) = \widehat \bu\left(\frac{\om}{\ep^2},Z\right) $
satisfies 
\begin{align}
\frac{\partial ^2 \widehat \bu_\ep(\om,Z)}{\partial Z^2} +
\frac{\om^2}{\ep^4c_o^2}\Big[ 1+ \ep \tilde\mu \Big(
  \frac{Z/c_o}{\ep^2} \Big) \Big] \widehat \bu_\ep(\om,Z) = 0,
\qquad Z > 0,
\end{align}
with initial conditions  
\begin{align}
\widehat \bu_\ep(\om,0) = \frac{\widehat f(\om)}{ 2 \sqrt{\zeta_o} },
\qquad \label{eq:RO} \frac{\partial \widehat \bu_\ep(\om,0)}{\partial
  Z} = -\frac{i \om \widehat f(\om)}{ 2 c_o\sqrt{\zeta_o}}. \nonumber
\end{align} 
This Cauchy problem is studied in \cite[Section 7.5.3]{newbible2007},
for random harmonic oscillators. The rate corresponding to exponential
growth in $Z$ of the square root of the energy $|\hat u_\ep|^2 +
|\partial_Z \hat u_\ep|^2$ is the Lyapunov exponent defined in
\cite[Proposition 7.9]{newbible2007}. It coincides with the rate
$\om^2 |\gamma^{(c)}|$ predicted by our analysis.

\subsection{Example of a medium in translational motion}
\label{sect:EX_R}
The integral equation (\ref{eq:inteq}) of the wave front is
complicated when $\alpha = \beta = 2$, because of the strong
dependence of the random travel time fluctuations on the starting
point $s$ of the characteristics. We do not have an explicit
description of the wave front for an arbitrary model $\mu$ of the
fluctuations, but there are models where the problem simplifies.  We
illustrate here one of them, corresponding to media in uniform
translational motion
\begin{equation}
\mu\Big(\frac{t}{\ep^2},\frac{z}{\ep^2}\Big) = \eta \Big(\frac{z-v
t}{\ep^2}\Big),
\label{eq:EX1.1}
\end{equation}
where $\eta$ is a mean zero stationary random process  with
autocorrelation $\varPhi$. The process $\eta$ is as smooth as before and with
finite dependence range.   The speed $v$ may be positive or negative, but
it satisfies the inequality $|v| < c_o$.

The random travel time fluctuations (\ref{eq:Rshift}) are
\begin{eqnarray}
\cW_\ep(\tau,s) &=& \frac{1}{\ep} \int_0^\tau du \, \eta
\Big(\frac{(c_o-v)u}{\ep^2} - v s - v \cW_\ep(u,s) \Big) ,
\label{eq:EX1.5}
\end{eqnarray}
and changing variables $ \zeta = (c_o-v) u - \ep^2 v s - \ep^2 v
\cW_\ep(u,s), $ with $\zeta$ lying in the interval $ \zeta \in \left[
  O(\ep^2), (c_o-v) \tau + O(\ep^2)\right]$, we obtain
\begin{equation}
\cW_\ep(\tau,s) = \frac{1}{\ep (c_o-v)} \int_0^{(c_o-v)\tau} d \zeta \, 
\eta \left({\zeta}/{\ep^2} \right) +  O(\ep).
\label{eq:EX1.6}
\end{equation}
This is independent of $s$ to leading order, and $\cW_\ep$ converges
in distribution, as $\ep \to 0$, to the diffusion process $ \cW(\tau)
= \cD \cB(\tau), $ where $\cB(\tau)$ is standard Brownian motion and $
\cD =\sqrt{{\widehat \varPhi(0)}/{(c_o-v)}}.  $ Thus, we can set
$Y_{\alpha,\beta}^\ep = 0$ and $\cS_\ep(y;\tau,s) = y$ in
(\ref{eq:inteq}), and the random processes (\ref{eq:nupm}) are
\begin{equation}
\mu_{\alpha,\beta}^\pm(t,z) = (1 \mp {v}/{c_o}) \eta'(z-vt).
\end{equation}
The pulse stabilization results in  Theorem  \ref{thm.1} and Corollary
\ref{cor.1s} hold as stated, with
\begin{equation}
\Psi(s) = \Big[(v/c_o)^2-1\Big] \varPhi''\left((c_o+v)s\right),
\end{equation}
and the parameter $\gamma(\om)$ in the pulse shaping kernel is 
\begin{align}
\gamma(\om) &= \frac{(v^2-c_o^2)}{4 \om^2 c_o} \int_0^\infty du \,
\varPhi''\left((c_o+v)u\right) e^{2i \om u} \nonumber \\ &= i
\frac{(v-c_o)}{2 \om c_o (c_o + v)} + \gamma^{(c)}(\om) + i \gamma^{(s)}(\om).
\end{align}
Here we used integration by parts and the normalization $\varPhi(0) = 1$, and let 
\begin{align}
\gamma^{(c)}(\om) &= \frac{(c_o-v)}{2 c_o (c_0+v)^2} \hat
\varPhi\left(\frac{2 \om}{c_o +v}\right) \\ 
\gamma^{(s)}(\om) &= { \frac{(c_o-v)}{ c_o (c_0+v)^2} } \int_0^\infty
dz \,\varPhi(z) \sin\left[2 \om z/(c_o+v) \right].
\label{eq:gammaseta}
\end{align}

These results show that as $v \uparrow c_o$ the transformation of the
pulse shape becomes small, as if the pulse had propagated only through
a thin section of the random medium. This is because the right
propagating wave component travels with a speed only slightly larger
than the medium.  However, the random travel time correction becomes
large since the averaging becomes ``less effective''.  The
diffusion becomes  small also in the limit 
$v \downarrow-c_o$, due to the assumed rapid decay of the power spectral density of the
fluctuations, while the dispersive effect becomes strong.

\section{Proof of results}
\label{sect:proof}
The proof builds on the approach in \cite[Section
  8.1]{newbible2007}. See also
\cite{burridge1989multimode,burridge1988one}. We use random travel
time coordinates to derive in section \ref{sect:proof.inteq} a time
domain integral equation for the wave front. The coordinates are
motivated by the method of characteristics applied to the first order
system of partial differential equations satisfied by the right and
left going waves, as described in section \ref{sect:proof.char}.  The
integral equation for the wave front is (\ref{eq:inteq}) and holds for
all $\alpha, \beta \in (0,2]$.  However, it simplifies in regimes 1-2,
where the random fluctuations of the travel time are independent of
the characteristics in the limit $\ep \to 0$. The simplification is
explained below and it leads to the explicit
expression of the wave front stated in Theorem \ref{thm.1}, as shown
in section \ref{sect:slow}.
\subsection{The wave decomposition and random travel time coordinates}
\label{sect:proof.char}
The right and left going waves $A_\ep(t,z)$ and $B_\ep(t,z)$ defined
by (\ref{eq:dec.1}) satisfy the following system of first order
partial differential equations derived in appendix \ref{ap:A}
\begin{eqnarray}
\frac{1}{c_\ep(t,z)} \frac{\partial A_\ep(t,z)}{\partial t} +
\frac{\partial A_\ep(t,z)}{\partial z} &=& \frac{1}{\ep} \left[
M_\ep(t,z) + N_\ep(t,z) \right] B_\ep(t,z), \label{eq:char.1} \\
\frac{1}{c_\ep(t,z)} \frac{\partial B_\ep(t,z)}{\partial t} -
\frac{\partial B_\ep(t,z)}{\partial z} &=& \frac{1}{\ep} \left[
-M_\ep(t,z) + N_\ep(t,z) \right] A_\ep(t,z), \label{eq:char.2}
\end{eqnarray}
with initial condition (\ref{eq:dec.2}) for $A_\ep$ and final
condition (\ref{eq:dec.3}) for $B_\ep$. Here 
\begin{equation}
M_\ep(t,z) = -
\frac{\ep^{2-\beta}\mu_z\left(\frac{t}{\ep^\alpha},\frac{z}{\ep^\beta}\right)}{
  2 \left[1 + \ep \mu
    \left(\frac{t}{\ep^\alpha},\frac{z}{\ep^\beta}\right)\right]},
\label{eq:char.3}
\qquad 
N_\ep(t,z) = -\frac{
\ep^{2-\alpha}\mu_t\left(\frac{t}{\ep^\alpha},\frac{z}{\ep^\beta}\right)}{
2 c_o \left[1 + \ep \mu
\left(\frac{t}{\ep^\alpha},\frac{z}{\ep^\beta}\right)\right]^{1/2}},
\end{equation}
where $\mu_z$ and $\mu_t$ are the partial derivatives of $\mu$ in $z$
and $t$.

We solve the first order system (\ref{eq:char.1})-(\ref{eq:char.2})
using the method of characteristics.  The characteristics of the right
going waves are the curves $\left(t_\ep(z,s),z\right)$ in the $(t,z)$
plane, where $t_\ep(z,s)$ is the travel time (\ref{eq:tt}), and $s$
parameterizes the impinging pulse at $z = 0$, giving the initial
condition
\begin{equation}
\label{eq:iniAnew}
A_\ep(t = \ep^2 s,z=0) = f(s), \qquad s \in [0,S].
\end{equation}  
The right going wave observed at the random travel time is obtained by
integrating equation (\ref{eq:char.1}) along the characteristics
\begin{equation}
A_\ep(t_\ep(z,s),z) = f(s) + \frac{1}{\ep} \int_0^z d z' \,
  G_\ep(t_\ep(z',s),z'),
\label{eq:char.6new}
\end{equation}
where
\begin{equation}
G_\ep(t,z) = \left[M_\ep(t,z) + N_\ep(t,z)\right] B_\ep(t,z).
\label{eq:char.7new}
\end{equation}

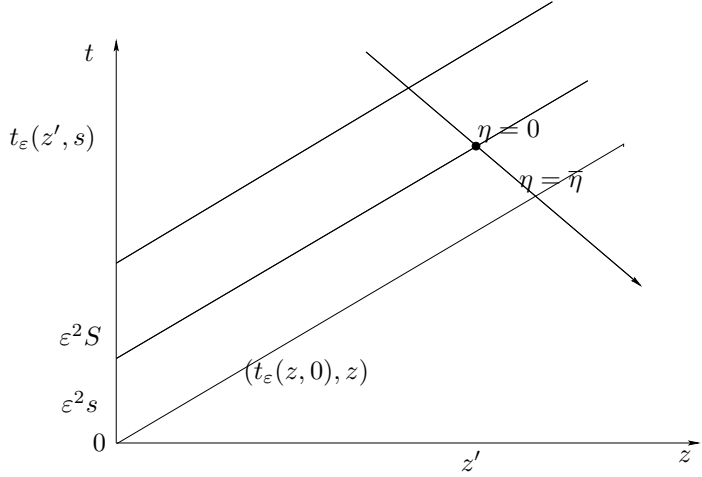
\begin{figure}[t]
\vspace{-0.4in}
  \begin{picture}(0,0)%
\raisebox{0.in}{    
\centerline{\includegraphics[width=0.6\textwidth]{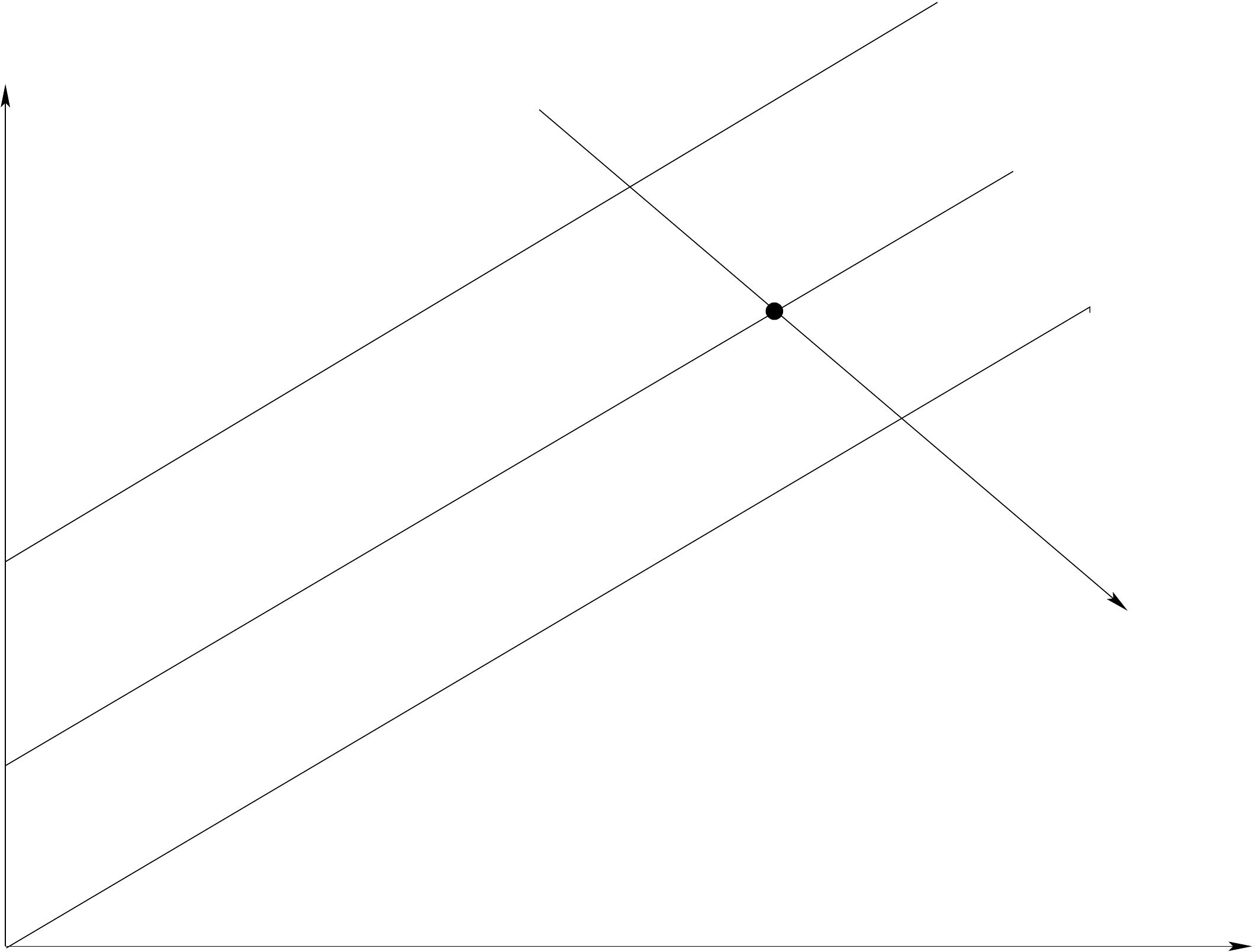}}}%
\end{picture}%
\setlength{\unitlength}{4144sp}%
\begingroup\makeatletter\ifx\SetFigFontNFSS\undefined%
\gdef\SetFigFontNFSS#1#2#3#4#5{%
  \reset@font\fontsize{#1}{#2pt}%
  \fontfamily{#3}\fontseries{#4}\fontshape{#5}%
  \selectfont}%
\fi\endgroup%
\begin{picture}(2500,2900)(0,400)
  \put(1050,2700){\makebox(0,0)[lb]{\smash{{\SetFigFontNFSS{12}{14.4}{\rmdefault}{\mddefault}{\updefault}{\color[rgb]{0,0,0}{\normalsize
              ${t}$} }%
        }}}}
  \put(900,1000){\makebox(0,0)[lb]{\smash{{\SetFigFontNFSS{12}{14.4}{\rmdefault}{\mddefault}{\updefault}{\color[rgb]{0,0,0}{\normalsize
              ${\ep^2 S}$} }%
        }}}}
  \put(600,2200){\makebox(0,0)[lb]{\smash{{\SetFigFontNFSS{12}{14.4}{\rmdefault}{\mddefault}{\updefault}{\color[rgb]{0,0,0}{\normalsize
              $t_\ep(z',s)$} }%
        }}}}
  \put(850,600){\makebox(0,0)[lb]{\smash{{\SetFigFontNFSS{12}{14.4}{\rmdefault}{\mddefault}{\updefault}{\color[rgb]{0,0,0}{
              \normalsize $\ep^2 s$} }%
        }}}}
  \put(1100,370){\makebox(0,0)[lb]{\smash{{\SetFigFontNFSS{12}{14.4}{\rmdefault}{\mddefault}{\updefault}{\color[rgb]{0,0,0}{\normalsize
              $0$} }%
        }}}}
  \put(4600,300){\makebox(0,0)[lb]{\smash{{\SetFigFontNFSS{12}{14.4}{\rmdefault}{\mddefault}{\updefault}{\color[rgb]{0,0,0}{\normalsize
              $z$} }%
        }}}}
  \put(3300,250){\makebox(0,0)[lb]{\smash{{\SetFigFontNFSS{12}{14.4}{\rmdefault}{\mddefault}{\updefault}{\color[rgb]{0,0,0}{\normalsize
              $z'$} }%
        }}}}
  \put(3400,2250){\makebox(0,0)[lb]{\smash{{\SetFigFontNFSS{12}{14.4}{\rmdefault}{\mddefault}{\updefault}{\color[rgb]{0,0,0}{\normalsize
              $\eta = 0$} }%
        }}}}
  \put(3650,1950){\makebox(0,0)[lb]{\smash{{\SetFigFontNFSS{12}{14.4}{\rmdefault}{\mddefault}{\updefault}{\color[rgb]{0,0,0}{\normalsize
             $\eta = \overline{\eta}$} }%
        }}}}
  \put(2000,800){\makebox(0,0)[lb]{\smash{{\SetFigFontNFSS{12}{14.4}{\rmdefault}{\mddefault}{\updefault}{\color[rgb]{0,0,0}{\normalsize
              $\left(t_\ep(z,0),z\right)$} }%
        }}}}
\end{picture}%
\vspace{0.05in}
\caption{The characteristics curves $\left(t_\ep(z,s),z\right)$ of the
  right going waves start at times $\ep^2 s$ on the initial curve $z =
  0$, with $s \in [0,S]$ parameterizing the pulse. We draw the
  characteristics approximately, as straight lines of slope $1/c_o$, but
  we recall that $t_\ep(z,s)$ has random fluctuations of order
  $\ep^2$. The characteristic of the left going wave passing through
  the point $\left(t_\ep(z',s),z'\right)$ is parameterized by $-\eta$,
  with $\eta = 0$ at the point of intersection. The front of the pulse
  is on the characteristic curve $\left(t_\ep(z,0),z\right)$, stemming
  from $s = 0$ at $z = 0$. }
\label{fig:1}
\end{figure}

To determine the left going wave $B_\ep\left(t_\ep(z',s),z'\right)$,
we integrate (\ref{eq:char.2}) along its characteristic passing
through $\left(t_\ep(z',s),z'\right)$. As illustrated in
Figure \ref{fig:1}, we parameterize this characteristic by $-\eta$,
starting from the point of intersection. We have
\begin{align}
  B_\ep(t_{\ep}(z',s)+\tilde{t}_\ep(\eta; z',s), z'+\eta) = &
  B_\ep(t_{\ep}(z',s), z') - \nonumber \\ & \frac{1}{\ep} \int_0^\eta
  d \eta' \, Q_\ep( t_{\ep}(z',s) + \tilde{t}_{\ep}(\eta'; z',s),
  z'+\eta') ,
\label{eq:char.12new}
\end{align}
where
\begin{equation}
Q_\ep(t,z) = \left[-M_\ep(t,z) + N_\ep(t,z)\right] A_\ep(t,z),
\label{eq:char.10new}
\end{equation}
and 
\begin{equation}
\label{eq:char.10new1}
  \tilde{t}_{\ep}(\eta; z',s)   = - \int_0^{\eta} 
  \frac{d\eta'}{c_\ep \left(t_{\ep}(z',s)+ 
      \tilde{t}_{\ep}(\eta'; z',s), z'+\eta' \right)} .
\end{equation}
By causality, the left going wave is zero ahead of the front given by
the curve $\left(t_\ep(z,0),z\right)$ for $z>0$.  We denote 
$\eta$ at the front by $\overline{\eta}$.  It is defined implicitly by
\begin{align*}
\tilde{t}_\ep(\overline{\eta};z',s) &= t_\ep(z'+\overline{\eta},0) -
t_\ep(z',s) = \frac{\overline{\eta}}{c_o} - \ep^2 \Big[ s + \cW_\ep
  \Big(\frac{z' }{c_o},s\Big) - \cW_\ep \Big(\frac{z'+
    \overline{\eta}}{c_o},0\Big)\Big].
\end{align*}
Now substitute this in (\ref{eq:char.10new1}) and use the model
(\ref{eq:as.5}) of the wave speed to obtain
\begin{align*}
\frac{2 \overline{\eta}}{c_o} + \ep \int_0^{\overline{\eta}} \frac{d
  \eta'}{c_o} \mu \Big( \frac{t_\ep(z',s) +
  \tilde{t}_\ep(\eta';z',s)}{\ep^\alpha},
\frac{z'+\eta'}{\ep^\beta}\Big) = \ep^2 s + \\ \ep^2 \Big[\cW_\ep
  \Big(\frac{z'}{c_o},s\Big) - \cW_\ep
  \Big(\frac{z'+\overline{\eta}}{c_o},0\Big) \Big].
\end{align*}
Recalling that $\mu$ is bounded, and using definition
(\ref{eq:Rshift}) of $\cW_\ep$ and the mean value theorem we have 
\[
\overline{\eta} \Big[ \frac{2}{c_o} + O(\ep)\Big] = \ep^2 s + \ep^2
\Big[\cW_\ep \Big(\frac{z'}{c_o},s\Big) - \cW_\ep
  \Big(\frac{z'}{c_o},0\Big) \Big].
\]
Here $O(\ep)$ denotes a term that is bounded  by $C \ep$,
for an $\ep$ independent constant $C$. Thus, we conclude that 
\begin{equation}
  \overline{\eta} = \frac{\ep^2 c_o}{2} \left[ s +
    \cW_\ep\left(\frac{z'}{c_o},s\right)-
    \cW_\ep\left(\frac{z'}{c_o},0\right)\right] + O(\ep^3).
\label{eq:overlineeta}
\end{equation}
Lemma \ref{thm.01} states that in regimes 1-2 the random process
$\cW_\ep(\tau,s)$ is independent of $s$ to leading order, so
$\overline{\eta} \approx \ep^2 c_o s/2$.  In the general case $\alpha,
\beta \in (0,2]$ we may have $\cW_\ep(\tau,s)$ depending on $s$, so
$\overline{\eta}$ may be random.  Nevertheless, $\overline{\eta}$ is
guaranteed to be positive to leading order, because the characteristic
stemming from $\ep^2 s$ at $z=0$ cannot intersect\footnote{The
  characteristics do not intersect because they are solutions of a
  flow (ordinary differential equation) problem that is uniquely
  solvable \cite[Section 3.2.5.a]{Evans}.}  the front of the pulse,
the characteristic stemming from the origin.

The left going wave at the point $(t_\ep(z',s),z')$ of intersection
with the right going characteristic is 
\begin{equation}
  B_\ep(t_{\ep}(z',s), z')  =
  \frac{1}{\ep}
  \int_0^{\overline{\eta}} d \eta' \, Q_\ep(
  t_{\ep}(z',s) + \tilde{t}_{\ep}(\eta'; z',s), z'+\eta').
\label{eq:char.12new1}
\end{equation}
It follows by setting $\eta = \overline{\eta}$ in (\ref{eq:char.12new}),
because the left hand side is the left going wave which vanishes at
the wave front.

\subsection{The integral equation}
\label{sect:proof.inteq}
The integral equation for the wave front is obtained by substituting
(\ref{eq:char.12new1}) in (\ref{eq:char.6new}). We write it in more
convenient form by changing variables in (\ref{eq:char.12new1}) as
$\eta' = \ep^2 c_o(s-y)/{2}$, with $y$ satisfying 
\begin{equation}
Y^{\ep}_{\alpha,\beta}(z',s) \le y \le s, \quad \mbox{for}~ ~ Y^\ep_{\alpha,\beta}(z',s) =
\cW_\ep\left({z'}/{c_o},0\right) - \cW_\ep\left({z'}/{c_o},s\right).
\label{eq:char.12new4}
\end{equation}

The equation for the left going wave becomes
\begin{align}
  \hspace{-0.1in}B_\ep(t_{\ep}(z',s), z') = \frac{\ep c_o}{2}\hspace{-0.05in}
  \int_{Y^\ep_{\alpha,\beta}(z',s)}^{s} \hspace{-0.08in}d y \, Q_\ep\Big[
    t_{\ep}(z',s) + \tilde{t}_{\ep}\Big(\frac{\ep^2 c_o(s-y)}{2};
    z',s\Big), z'+\frac{\ep^2 c_o(s-y)}{2}\Big],
\label{eq:DefB}
\end{align}
and substituting it in (\ref{eq:char.6new}) and using definition
(\ref{eq:char.10new}) of $Q_\ep$, we get the following equation for
the right going wave
\begin{align}
  A_\ep\left[t_\ep(z,s),z\right] &= f(s) - \frac{c_o}{2} \int_0^z dz'
  \left( M_\ep + N_\ep \right) \left( t_\ep(z',s), z'\right) \times
  \nonumber \\ & \hspace{0.in}\int_{Y_{\alpha,\beta}^\ep(z',s)}^{s} d y \,
  A_\ep\Big[ t_\ep(z',s) + \tilde{t}_{\ep}\Big(\frac{\ep^2
      c_o(s-y)}{2}; z',s\Big), z'+\frac{\ep^2 c_o(s-y)}{2}
    \Big]\times \nonumber \\ & \hspace{0.in}\Big(M_\ep - N_\ep
  \Big) \Big[ t_\ep(z',s) + \tilde{t}_{\ep}\Big(\frac{\ep^2
      c_o(s-y)}{2}; z',s\Big), z'+\frac{\ep^2 c_o(s-y)}{2}\Big] .
     \label{eq:intnew1}
\end{align}
Next we rewrite this integral equation in terms of the right propagating wave in (\ref{eq:wf}) expressed in 
characteristic coordinates and to do so we introduce below a coordinate mapping. 
First we  obtain from (\ref{eq:char.10new1}) that
\begin{align}
  \tilde{t}_{\ep}\left({\ep^2 c_o(s-y)}/{2}; z',s\right) &= -
  \int_0^{\frac{\ep^2 c_o(s-y)}{2}} \frac{d\eta'}{c_\ep
    \left(t_{\ep}(z',s)+ \tilde{t}_{\ep}(\eta'; z',s), z'+\eta'
    \right)} \nonumber \\
 &= -{\ep^2(s-y)}/{2} + \ep^3 R_\ep(y,z'/c_o,s),
\label{eq:defBp}
\end{align}
where we used definition (\ref{eq:as.5}) of $c_\ep$, changed variables
$\eta' = \ep^2 \zeta$, and introduced the remainder
\begin{equation*}
R_\ep(y,z'/c_o,s) = -\int_0^{\frac{c_o(s-y)}{2}} d \zeta \, \mu \Big(
\frac{t_\ep(z',s) + \widetilde t_\ep(\ep^2 \zeta;z',s)}{\ep^\alpha},\frac{z'+
\ep^2 \zeta}{\ep^\beta}\Big),
\end{equation*}
with both $R_\ep$ and $\partial_y R_\ep$ bounded.
Recalling definition (\ref{eq:tt}) of $t_\ep$, we get 
\begin{align}
t_\ep(z',s) + \tilde{t}_{\ep}\Big(\frac{\ep^2 c_o(s-y)}{2}; z',s\Big)
=& \frac{z'}{c_o} + \frac{\ep^2(s-y)}{2} + \nonumber \\ &\ep^2 \Big[y +
  \cW_\ep\Big(\frac{z'}{c_o},s\Big) + \ep R_\ep(y,z'/c_o,s)\Big],
\label{eq:Ad1}
\end{align}
and rewrite the right hand side in terms of the mapping
$\cS_\ep(y;\tau,s)$ defined by
\begin{align}
  \cS_\ep\Big(y;\frac{z'}{c_o},s\Big) &= y +
  \cW_\ep\Big(\frac{z'}{c_o},s\Big) - \cW_\ep\Big(\frac{z'}{c_o}+
  \frac{\ep^2 (s-y)}{2} ,\cS_\ep\Big(y;\frac{z'}{c_o}, s\Big) \Big) +
  \ep R_\ep(y,z'/c_0,s).
\label{eq:Ad2}
\end{align}
 Since $\ep^{2}(s-y)\partial_{\tau} \cW_\ep(\tau,s) = O(\ep)$,
  uniformly in $\tau, s$ and $y$, we see that this is the same
as definition (\ref{eq:mapSe}), with $\cR_\ep$ equal to $R_\ep$ plus
an order one correction coming from expanding $\cW_\ep$ around $\tau'
= z'/c_o$.   Using (\ref{eq:Ad2}) in (\ref{eq:Ad1}) we obtain that 
\begin{align*}
A_\ep\Big[ t_\ep(z',s) + \tilde{t}_{\ep}\Big(\frac{\ep^2 c_o(s-y)}{2};
  z',s\Big), z'+\frac{\ep^2 c_o(s-y)}{2} \Big] =&
\\ & \hspace{-0.7in}a_\ep\Big(z'/c_o +
\frac{\ep^2(s-y)}{2},\cS_\ep\left(y;z'/c_o,s\right) \Big),
\end{align*}
and equation (\ref{eq:intnew1}) becomes 
\begin{align}
 a_\ep(\tau,s) \nonumber &= f(s) - \frac{c_o^2}{2} \int_0^\tau d\tau'
 \left( M_\ep + N_\ep \right) \left( t_\ep(c_o \tau',s), c_o
 \tau'\right) \times \nonumber
 \\ & \hspace{0.1in}\int_{Y_{\alpha,\beta}^\ep(c_0 \tau',s)}^{s} d y \,
 a_\ep\left(\tau' + {\ep^2(s-y)}/{2},\cS_\ep\left(y;\tau',s\right)
 \right)\times \nonumber \\ & \hspace{0.1in}\left(M_\ep - N_\ep
 \right) \Big[ 
 t_\ep(c_o \tau',s) + \tilde{t}_{\ep}\Big(\frac{\ep^2
      c_o(s-y)}{2}; z',s\Big)
 , c_o
   \tau'+\frac{\ep^2 c_o(s-y)}{2}\Big]  .
 \label{eq:intnew2}
\end{align}
Equation (\ref{eq:inteq}) follows by substitution of
definitions (\ref{eq:char.3}) of $M_\ep$ and $N_\ep$.

It remains to show that (\ref{eq:intnew2}) simplifies in regimes $1$ and $2$. 
Note first that in view of  Lemma \ref{thm.01}, $Y_{\alpha,\beta}^\ep$ is smaller 
than $\ep$ to some positive power, on  an  event  $E$ with  probability 
$\PP(E) > 1 - O(\ep^\theta)$.  Recall  also that the  bound $Y_{\alpha,\beta}^\ep$ 
on $y$ corresponds to the wavefront, meaning that for $y < Y_{\alpha,\beta}^\ep$ the contribution to the 
integral is zero. Thus, we can extend the y integral to the interval $[-S,s]$ without making any mistake. 
In view of Lemma \ref{thm.01}  and Eq. (\ref{eq:Ad2}) we also have that:
 \begin{equation}
\PP \Big( \sup_{\tau \in [0,T], 0 \le y \le s \le S} |\partial_y
\cS^\ep(y;\tau,s) - 1 | \ge c\ep   \Big) \le C \ep^{\theta}. 
\label{eq:cor2e}
\end{equation}
Then we get via a change of variables in (\ref{eq:intnew2}) that on an event with probability  
 $1 - O(\ep^\theta)$ , 
\begin{align}
 a_\ep(\tau,s) \nonumber &= f(s) - \frac{c_o^2}{2} \int_0^\tau d\tau'
 \left( M_\ep + N_\ep \right) \left( t_\ep(c_o \tau',s), c_o
 \tau'\right) \times \nonumber
 \\ & \hspace{0.1in} \int_{-S}^{s} d y \,
 a_\ep\left(\tau' + {\ep^2(s-y)}/{2},y
 \right)\times \nonumber \\ & \hspace{0.1in}\left(M_\ep - N_\ep
 \right) \Big[ 
 t_\ep(c_o \tau',s) + \tilde{t}_{\ep}\Big(\frac{\ep^2
      c_o(s-y)}{2}; z',s\Big)
 , c_o
   \tau'+\frac{\ep^2 c_o(s-y)}{2}\Big]    \times \nonumber
 \\ & \hspace{0.1in}
   (1 +  \ep   h_\ep(y,\tau',s)) ,
 \label{eq:intnew3}
\end{align}
with some bounded $h_\ep$. 
 
 We show next that $\partial_\tau a_\ep$ is bounded, so we can simplify the 
 first argument of $a_\ep$ in (\ref{eq:intnew3}).
When taking the partial derivative with respect to $\tau$  we see that we can bound $\partial_\tau a_\ep$ in
terms of $a_\ep$. To bound the latter, we recall that $M_\ep$ and
$N_\ep$ are uniformly bounded.   
 We then  get 
 \[
\|a_\ep(\cdot,s)\|_{\infty}   \le \|f\|_\infty +
C \int_{-S}^s  dy \|a_\ep(\cdot,y)\|_{\infty},
\]
where we introduced the notation 
\[
\|a_\ep (\cdot,s)\|_{\infty} = \sup_{\tau} |a_\ep(\tau,s)|.
\]
The uniform boundedness of $a_\ep$, and
therefore of $\partial_\tau a_\ep$, follows from Gronwall's
inequality.
We thus find, using also the smoothness of $\mu$,  that   on an event with probability  
 $1 - O(\ep^\theta)$  
   we can simplify
the equation (\ref{eq:inteq}) of the wave front as
\begin{align}
a_\ep(\tau,s) = f(s) + \int_0^\tau \hspace{-0.05in} d
\tau' \hspace{-0.05in}\int_{-S}^s \hspace{-0.05in} dy \,
\mathcal{G}_\ep(\tau',s,y) a_\ep(\tau',y) + O(\ep),
\label{eq:SIMPL1}
\end{align}
where 
\begin{align} 
\mathcal{G}_\ep(\tau',s,y) := & -\frac{c_o^2}{8}\, \mu_{\alpha,\beta}^+
\hspace{-0.03in}\Big[\frac{\tau' + \ep^2\left(s +
    \cW_\ep(\tau',s)\right)}{\ep^\alpha}, \frac{c_o \tau'}{\ep^\beta}
  \Big] \times \nonumber \\ &
\mu_{\alpha,\beta}^- \hspace{-0.03in}\Big[\frac{ \tau' + \ep^2
    \left(\frac{y+s}{2}+ \cW_\ep(\tau',s)\right)}{\ep^\alpha},
  \frac{c_o \tau' + \ep^2 c_o
    (s-y)/2}{\ep^\beta}\Big], \label{eq:SIMPL2}
\end{align}
with $\mu_{\alpha,\beta}^\pm$ defined in (\ref{eq:nupm}).

\subsection{Stochastic averaging}
\label{sect:slow}
To show that $a_\ep \to \overline{a}$ as stated in Theorem 
\ref{thm.1}, it suffices to show that 
\begin{equation}
\label{eq:pf.2}
    \sup_{\tau \le T}  \| a_\ep(\tau,\cdot) -
\overline{a}(\tau,\cdot) \|_{\infty}   \stackrel{\ep \to
  0}{\longrightarrow} 0.
\end{equation}
in probability. 
We give the proof in regime 1, with $\alpha < 1$ and $\beta = 2$, but
the analysis in regime 2 is very similar.
\subsubsection{The linear integral operators}
Let us rewrite equation (\ref{eq:SIMPL1}) as
\begin{equation}
a_\ep(\tau,s) =  f(s) +
\int_0^\tau du \, \left[ \cLe(u) a_\ep(u,\cdot) \right] (s),
\label{eq:pf.11}
\end{equation}
using the linear operator $\cLe(u)$ defined by 
\begin{align}
  \left[\cLe(u) \varphi\right](s) &= \int_{-S}^s dy\,
  \mathcal{G}_\ep(u,s,y) \varphi(y), \qquad \forall \, \varphi \in
  \cC_{[-S,S]},
\label{eq:pf.4}
\end{align}
on the space $\cC_{[-S,S]}$ of continuous functions on $ [-S,S]$.
 Here $\mathcal{G}_\ep$ is defined by
(\ref{eq:SIMPL2}) for $\tau' = u$, and $\mu_{\alpha,2}^\pm$ is
obtained from definition (\ref{eq:nupm})
\begin{equation}
\mu^\pm_{\alpha,2}(t,z) = \mu_z(t,z) \pm
\ep^{2-\alpha}c_o^{-1}\mu_t(t,z).
\label{eq:nupm2}
\end{equation}
Since $\mu_{z}$ is bounded, we conclude that $\cLe(u)$ is Lipschitz
continuous (i.e., bounded) with deterministic constant
$C_{_\mathcal{L}}$ that is independent of $\ep$ and $u$,
\begin{equation}
\|\cLe(u) \varphi \|_{\infty} \le C_{_\mathcal{L}} \|\varphi
\|_{\infty}, \qquad \forall \varphi \in \cC_{[-S,S]}.
\label{lemma.pf.1}
\end{equation}
Similarly, we write the equation satisfied by the limit $\overline{a}$
as
\begin{equation}
\overline{a}(\tau,s) = f(s) + \int_{0}^\tau du \, \left[ \, \cLo \, 
\overline{a}(u,\cdot) \right] (s),
\label{eq:pf.12}
\end{equation}
using the bounded linear operator $\cLo:\cC_{[-S,S]} \to \cC_{[-S,S]}$,
defined by
\begin{equation}
\left[\cLo\, \varphi\right](s) =-\frac{c_o^2}{8} \int_{-S}^{s} d y \,
\varphi(y) \, \Psi \left(\frac{s-y}{2}\right), \qquad \forall \, \varphi \in
  \cC_{[-S,S]},
\label{eq:pf.5}
\end{equation}
where
\begin{equation}
\Psi( s ) = \EE \left\{ \mu_z(t,z) \mu_z({t,z + c_o s})\right\} =
-\partial^2_z \Phi(0,c_o s).
\label{eq:pf.6}
\end{equation}
Note that (\ref{eq:pf.5}) is the same as the expectation of
(\ref{eq:pf.4}) up to a term of order $\ep^{2-\alpha}$.
We also remark that we could have replaced the lower integration bound 
$-S$ in (\ref{eq:pf.5}) by 0 due to the support of the initial condition $f$  in
(\ref{eq:pf.12}).

\subsubsection{Averaging}
\label{sect:proof.avg1}
We now show that the limit (\ref{eq:pf.2}) holds.
We have
\begin{equation} 
\label{eq:PP1}
a_\ep(\tau,s) -\overline{a}(\tau,s) = g_\ep(\tau,s) + 
\int_0^\tau du \left[ \cLe(u) \left(
 a_\ep(u,\cdot)-\overline{a}(u,\cdot)\right) \right](s),
\end{equation} 
where
\begin{align}
\label{eq:residual}
g_\ep(\tau,s) = \int_0^\tau du \left[ (\cLe(u) - \cLo
  )\overline{a}(u,\cdot)\right](s).
\end{align}
Recalling that $\cLe(u)$ is bounded, and using the triangle
inequality, we get from (\ref{eq:PP1})
\begin{equation}
\|a_\ep(\tau,\cdot) -\overline{a}(\tau,\cdot)\|_\infty \le \|
g_\ep(\tau,\cdot) \|_\infty + C_{_{\mathcal{L}}} \int_0^\tau du \, \|
a_\ep(u,\cdot) -\overline{a}(u,\cdot) \|_\infty.
\end{equation}
Now Gronwall's inequality gives {
\begin{equation}
\|a_\ep(\tau,\cdot) -\overline{a}(\tau,\cdot)\|_\infty \le e^{\tau
  C_{_{\mathcal{L}}} }\, \sup_{0<t<\tau}\| g_\ep(t,\cdot) \|_\infty ,
\label{eq:pf15}
\end{equation}}
and we have 
\begin{equation}
\sup_{0<\tau<T} \|a_\ep(\tau,\cdot) -\overline{a}(\tau,\cdot)\|_\infty \le  
  C' \, \sup_{0<t<T}\| g_\ep(t,\cdot) \|_\infty ,
\label{eq:pf15*}
\end{equation}
and $  \sup_{0<t<T} \| g_\ep(t,\cdot) \|_\infty $ is estimated
using the next lemma, proved in Appendix \ref{ap:D}.

\begin{lemma}
\label{lem.apD}
Suppose that $\mu(t,z)$ satisfies the mixing assumption (\ref{eq:asMix2}) and similar for $\mu_z$. Then for  $0< T' < T$
\begin{equation*}
\PP \left[ \sup_{\tau \in [0,\ep^q T], \, s \in [-S,S]}    \left|
  \int_{T'}^{T'+\tau} du \, [(\cLe(u)-\cLo)\varphi](s) \right|  >  c\ep^{\frac{\theta-\vartheta}{2}}  \right] 
  < C \ep^{ \vartheta +q} ,
 \end{equation*}
for any deterministic continuous function $\varphi(s)$, arbitrary $ q > 0$ and  $\vartheta$ satisfying 
 $max[0, 2(\theta-1)] <  \vartheta < \theta$ . 
Here  
$\theta$ is as in Lemma \ref{thm.01}.
\end{lemma}

Let us introduce a uniform discretization of the time, in steps
$\Delta$, and write the residual (\ref{eq:residual}) as 
\begin{align}
g_\ep(\tau,s) &= \sum_{j=0}^{\lfloor \tau/\Delta\rfloor - 1} \int_{j
  \Delta}^{(j+1)\Delta} du \left[ (\cLe(u) - \cLo
  )\overline{a}(u,\cdot)\right](s) + \nonumber \\ &\int_{\Delta
  \lfloor \tau/\Delta\rfloor}^{\tau} du \left[ (\cLe(u) - \cLo
    )\overline{a}(u,\cdot)\right](s).
\label{eq:NP1}
\end{align}
The last term is bounded by $C \Delta$ uniformly in $\tau$ and $s$,
because both $\cLe$ and $\cLo$ are bounded linear operators and $\bar
a(\tau,s)$ is   bounded.   Moreover, (\ref{lemma.pf.1}) and (\ref{eq:pf.12}) give
for $j \Delta \le u \le (j+1) \Delta$ that
\[
\| \cLe(u) \overline{a}(u,\cdot) - \cLe(u) \overline{a}(j \Delta,
\cdot) \|_{\infty} \le C_{\mathcal{L}} \|\overline{a}(u,\cdot) -
\overline{a}(j \Delta, \cdot)\|_{\infty} \le C (u - j \Delta),
\]
where $C$ is a constant that is independent of $j$, $\Delta$ and
$\ep$.  A similar result holds for $\cLo$, so we get 
\begin{align}
\|g_\ep(\tau,\cdot)\|_{\infty} \le & \sum_{j = 0}^{\lfloor
  \tau/\Delta\rfloor - 1} \sup_{s \in [-S,S]} \left|\int_{j
  \Delta}^{(j+1)\Delta} du \left[ (\cLe(u) - \cLo )\overline{a}(j
  \Delta,\cdot)\right](s)\right| + \nonumber \\ &C \sum_{j =
  0}^{\lfloor \tau/\Delta\rfloor - 1} \int_{j \Delta}^{(j+1)\Delta} du
\, (u - j \Delta) + C' \Delta, \label{eq:NP2}
\end{align} 
with redefined $\ep, \Delta$ independent constants $C$ and $C'$.
The second
term is
\[
\sum_{j = 0}^{\lfloor \tau/\Delta\rfloor - 1} \int_{j
  \Delta}^{(j+1)\Delta} du \, (u - j \Delta) =
\frac{1}{2}\sum_{j=0}^{\lfloor \tau/\Delta\rfloor - 1} {\Delta^2} \le
     {\tau \Delta }/{2},
\]
and thus together with the last term vanishes in the limit $\Delta \to
0$.  
In view of Lemma \ref{lem.apD} we get by choosing  $\Delta=\ep^q$  and 
$q = (\theta-\vartheta)/4$ that 
 \[
\PP \Big[  \sup_{\tau \le T} \|g_\ep(\tau,\cdot)\|_{\infty} > 
 c\ep^{\frac{\theta-\vartheta}{4}}  \Big]   <   C \ep^{\vartheta}  ,
\]
which gives the result. 

\section{{Summary}}
\label{sect:summary}
We have considered wave propagation in one dimensional random media in
the situation when the wave speed varies randomly in both space and
time.  We derived a mathematical framework involving a random integral
equation that can be used to characterize the pulse in this case.  The
main result is that the pulse stabilization phenomenon known in time
independent random media extends to the time variable case. Pulse
stabilization means that the amplitude of the wave transmitted from
the source to range $z$ has a deterministic form, when observed in a
time window centered at $\tau + \delta \tau$, and of width comparable
to the support of the initial pulse. Here $\tau$ is the deterministic
travel time to range $z$, and $\delta \tau$ is a small random
correction, on the scale of the support of the initial pulse. It is
well known that in time independent media the transmitted pulse is
faded and broadened, because scattering in the medium produces delays
and also transfers energy to the incoherent waves.  We assume random
temporal changes in the medium, modeled by temporal fluctuations of
the wave speed.  We call the time fluctuations slow when the life span
of a spatial realization of the random wave speed is longer than the
time of traversal of a correlation length and the pulse width. The
fluctuations are rapid when these time scales are similar.  The
results show that the slow fluctuations have to leading order no
effect on the wave front.  The rapid fluctuations affect both the
random arrival time correction $\delta \tau$ and the deterministic
pulse deformation. In particular, there is a trade-off between
frequency dependent attenuation due to scattering in the medium, and
energy gain induced by the temporal fluctuations.

The analysis is based on a stochastic time domain integral equation
satisfied by the wave front, and the pulse stabilization follows from
stochastic averaging.  Furthermore, the energy gain induced by the
temporal fluctuations is explained by drawing an analogy to the random
harmonic oscillator problem.
 
\section*{Acknowledgments}
The work of L. Borcea was partially supported by the  ONR Grant N00014-14-1-0077. Support from the AFOSR Grant
FA9550-12-1-0117 is also gratefully acknowledged.  The work of
  K. S{\o}lna was partially supported by AFSOR grant \# FA
  9550-11-1-0176.

\appendix 
\section{Derivation of the first order system for the right and left going 
waves}
\label{ap:A}
Let us take the $z$ derivative in the first equation in
(\ref{eq:dec.1}) to obtain
\begin{equation}
\frac{\partial A_\ep}{\partial z} = \zeta_\ep^{-1/2} \frac{\partial
\bp_\ep}{\partial z} - \frac{\bp_\ep}{2 \zeta_\ep^{3/2}} \frac{\partial
\zeta_\ep}{\partial z} + \zeta_\ep^{1/2} \frac{\partial
\bu_\ep}{\partial z} + \frac{\bu_\ep}{2 \zeta_\ep^{1/2}} \frac{\partial
\zeta_\ep}{\partial z}.
\label{eq:A.1}
\end{equation}
Since $\bp_\ep$ and $\bu_\ep$ satisfy (\ref{eq:form.1}), we have for
  $z \neq 0$
\begin{equation}
\frac{\partial A_\ep}{\partial z} = -\rho \zeta_\ep^{-1/2}
\frac{\partial \bu_\ep}{\partial t} - \frac{\bp_\ep}{2 \zeta_\ep^{3/2}}
\frac{\partial \zeta_\ep}{\partial z} - \zeta_\ep^{1/2} K_\ep^{-1}
\frac{\partial \bp_\ep}{\partial t} + \frac{\bu_\ep}{2 \zeta_\ep^{1/2}}
\frac{\partial \zeta_\ep}{\partial z}.
\label{eq:A.2}
\end{equation}
Equations (\ref{eq:dec.1}) give that $ \bp_\ep = \zeta_\ep^{1/2}
(A_\ep-B_\ep)/2$ and $\bu_\ep = \zeta_\ep^{-1/2} ({A_\ep+B_\ep})/2$,
so
\begin{eqnarray}
\frac{\partial \bp_\ep}{\partial t} &=& \frac{\zeta_\ep^{1/2}}{2} \Big(
\frac{\partial A_\ep}{\partial t}-\frac{\partial B_\ep}{\partial
t}\Big) + \frac{1}{4 \zeta_\ep^{1/2}} \frac{\partial
\zeta_\ep}{\partial t}\Big(A_\ep-B_\ep\Big) \nonumber \\
\frac{\partial \bu_\ep}{\partial t} &=& \frac{\zeta_\ep^{-1/2}}{2} \Big(
\frac{\partial A_\ep}{\partial t}+\frac{\partial B_\ep}{\partial
t}\Big) - \frac{1}{4 \zeta_\ep^{3/2}} \frac{\partial
\zeta_\ep}{\partial t}\left(A_\ep+B_\ep\right).
\label{eq:A.4}
\end{eqnarray}
Substituting in (\ref{eq:A.2}) and using that $ {\rho}/{\zeta_\ep}
= {1}/{c_\ep}, $ we obtain
\begin{equation}
\frac{\partial A_\ep}{\partial z} = -\frac{1}{c_\ep} \frac{\partial
  A_\ep}{\partial t} + \left(\frac{1}{2 \zeta_\ep} \frac{\partial
  \zeta_\ep}{\partial z} + \frac{1}{2 c_\ep \zeta_\ep} \frac{\partial
  \zeta_\ep}{\partial t}\right) B_\ep.
\label{eq:A.5}
\end{equation}
Equation (\ref{eq:char.1}) follows with 
\begin{equation}
M_\ep(t,z) = \frac{\ep}{2 \zeta_\ep} \frac{\partial
\zeta_\ep}{\partial z}, \qquad 
  N_\ep(t,z) =\frac{\ep}{2 c_\ep \zeta_\ep} \frac{\partial
    \zeta_\ep}{\partial t}.
\end{equation}
Equation (\ref{eq:char.2}) follows similarly. $\Box$
 
\section{Mixing assumption}
\label{ap:MIXING}  
For completeness  we show here  that the mixing conditions \eqref{eq:asMix} and \eqref{eq:asMix2} are implied by 
a $\phi-$mixing assumption, see also \cite{KP,PV}. To state the assumption, let $(\Omega,\cF,P)$ be a  probability 
space and denote by $\cF_{z_1}^{z_2}$ a family of  $\sigma-$algebras contained in $\cF$, so that
$\cF_{z_2}^{z_3} \subset  \cF_{z_1}^{z_4}, 0 \leq z_1 \leq z_2 \leq z_3 \leq z_4 \leq \infty$.
The $\phi-$mixing assumption is
\begin{equation}
\label{eq:PhiMix}
\sup_{z \ge 0} \, \sup_{U \in \cF_{z+\Delta z}^\infty, V \in \cF_0^z} \big|P(U|V)-P(U)| \le \big(\Delta z\big)^{-d}, \qquad \Delta z > 0,
\end{equation}
 with  $d > 2$. 
 We assume that conditional probabilities relative to $\cF_0^z, 0\leq z\leq\infty$, have a regular version so that with probability one we have the representation
 \[
 \EE\big[ \cdot \mid \cF_0^z \big](\omega') = \int_\Omega \cdot \PP_z(d\omega \mid \omega')  .
 \]
 Assume that $X$ is a real valued, bounded in absolute value by $c$ and   $\cF_0^z$ measurable random variable. Then
 \begin{eqnarray}
 \big| \EE\big[ X \mid \cF_0^z \big](\omega') -   \EE\big[ X \big] \big| &=&
 \big| \int_\Omega X(\omega) \big[ \PP_z(d\omega \mid \omega') -\PP(d\omega) \big] \big| \\
  &=&  
     \big| \int_\Omega X(\omega)  \upsilon_z(d\omega \mid \omega') \big| .
\end{eqnarray}
Here,   $\upsilon_z(U \mid \omega') =  \PP_z(U \mid \omega') -\PP(U)$ is a signed measure and,
from (\ref{eq:PhiMix}) and  the Hahn decomposition theorem, its variation 
$|\upsilon_z|(U \mid \omega')$ satisfies
\[
    \sup_{U \in \cF_{z+\Delta z}^\infty}  |\upsilon_z|(U \mid \omega') \leq 2  \big(\Delta z\big)^{-d} .
\]
Then we consequently have 
  \begin{eqnarray}
 \big| \EE\big[ X \mid \cF_0^z \big](\omega') -   \EE\big[ X \big] \big| &\leq& c 2  \big(\Delta z\big)^{-d}.
 \end{eqnarray}
 
 The bounds (\ref{eq:asMix}), respectively  (\ref{eq:asMix2}),  follow accordingly by considering
 the filtration generated by $\mu(\cdot,z)$ and the random variable $\mu(t,z)$,
  respectively $\mu(t_1,z_1) \mu(t_2,z_2)$ (moreover their derivatives).

\section{Proof of Lemma \ref{thm.01}}
\label{sect:Pf}
We prove the lemma in regime 1, where $\alpha < 1$ and $\beta =
2$.  The results extend to regime 2 using a change of variables, as explained in section 
\ref{sect:PfReg2}. 
We begin in section \ref{sect:Pf1} with the estimation
of the residual
\begin{equation}
q_\ep(z,s) = \cW_\ep\Big({z}/{c_o},s\Big) -W_\ep(z/c_o),
\label{eq:C1}
\end{equation}
where $W_\ep$ is defined in (\ref{eq:defchi}). The estimation of
$\partial_s \cW_\ep(z/c_o,s)$ is in section \ref{sect:Pf2}.

We use henceforth the notation $c$ and $C$ for constants that are independent of $\ep$. 

\subsection{Proof of estimate (\ref{lem2.eq2})}
\label{sect:Pf1}
Let us recall definitions (\ref{eq:Rshift}) and (\ref{eq:defchi}) to
write the residual (\ref{eq:C1}) as
\begin{align}
q_\ep (z,s) &= \ep^{1-\alpha} \int_0^z du \, \nu_\ep(u,s) \left[s +
  W_\ep((u-\ep^p)/c_o) + q_\ep(u-\ep^p,s)\right]+ \nonumber \\ &
\ep^{1-\alpha} \int_0^z du \, \widetilde \nu_\ep(u,s)
\left[W_\ep(u/c_o) - W_\ep((u-\ep^p)/c_o)+ q_\ep(u,s) -
  q_\ep(u-\ep^p,s)\right],
\label{eq:C3}
\end{align}
for $p$ satisfying
\begin{equation}
\label{eq:P}
{\alpha < p < 1 + \alpha/{2}}.
\end{equation}
By the mean value theorem
\begin{align*}
\mu\Big(\frac{u/c_o}{\ep^\alpha } + \ep^{2-\alpha}
\xi_\ep(u,s),\frac{u}{\ep^2}\Big) - \mu\Big(\frac{u/c_o}{\ep^\alpha}
+\ep^{2-\alpha} \eta_\ep(u,s),\frac{u}{\ep^2}\Big) = \ep^{2-\alpha}
\times \\\left[\xi_\ep(u,s)-\eta_\ep(u,s)\right]\mu_t
\Big(\frac{u/c_o}{\ep^\alpha} + \ep^{2-\alpha}
\zeta_\ep(u,s),\frac{u}{\ep^2}\Big),
\end{align*}
with $\zeta_\ep(u,s)$ between $\xi_\ep(u,s)$ and $\eta_\ep(u,s)$. In
the first term in (\ref{eq:C3}) we have
\begin{equation}
\label{eq:defnue}
\nu_\ep(u,s)= c_o^{-1}
\mu_t\Big(\frac{u/c_o}{\ep^\alpha} + \ep^{2-\alpha}
\zeta_\ep(u,s),\frac{u}{\ep^2}\Big),
\end{equation}
and $\zeta_\ep(u,s)$ is between $\xi_\ep(u,s) = s + W_\ep((u-\ep^p)/c_o)
+ q_\ep(u-\ep^p,s)$ and $\eta_\ep = 0$. The expression of $\widetilde
\nu_\ep$ in the second term of (\ref{eq:C3}) is similar, except that
$\zeta_\ep(u,s)$ is between $\xi_\ep(u,s) = s + W_\ep(u/c_o) +
q_\ep(u,s)$ and $\eta_\ep(u,s) = s + W_\ep((u-\ep^p)/c_o) +
q_\ep(u-\ep^p,s)$.  Using that $q_\ep(u,s)$ is supported for positive
$u$, we write equation (\ref{eq:C3}) as
\begin{align}
q_\ep (z,s) &= \ep^{1-\alpha} \hspace{-0.05in}\int_0^z du \,
\kappa_\ep(u,z,s) q_\ep(u,s) + \cT_{\ep,1}(z,s) + \cT_{\ep,2}(z,s),
\label{eq:C5}
\end{align}
with the notation
\begin{align*}
\cT_{\ep,1}(z,s) &= \ep^{1-\alpha}\hspace{-0.05in}\int_0^z
du \, \big[\widetilde \nu_\ep(u,s) - \nu_\ep(u,s)\big] \big[W_\ep(u/c_o) -
  W_\ep((u-\ep^p)/c_o)\big],
\nonumber \\
\cT_{\ep,2}(z,s) &= \ep^{1-\alpha}\hspace{-0.05in}\int_0^z du \,
\nu_\ep(u,s) \big[ s + W_\ep(u/c_o)\big],
\end{align*}
and integral kernel 
\begin{align*}
\kappa_\ep(u,z,s) &= \widetilde \nu_\ep(u,s) + \left[\nu_\ep(u+\ep^p,s)-
 \widetilde \nu_\ep(u+\ep^p,s)\right]1_{[0,z-\ep^p]}(u),
\end{align*}
where $1_{[m,M]}(u)$ is the indicator function of the interval
$[m,M]$.  We can bound $q_\ep$ using Gronwall's inequality, once we
estimate all the terms in (\ref{eq:C5}).

  Definition (\ref{eq:defchi}) and the
boundedness of $\mu$ and $\mu_t$ imply that
\begin{align*}
  \sup_{0 \le u \le z \le L} \ep^{1-\alpha} |\kappa_\ep(u,z,s)| \le C
  \ep^{1-\alpha}, \qquad \sup_{z\in [0, L]} 
  \big|\cT_{\ep,1}(z,s)\big| \le C \ep^{p -\alpha},
\end{align*}
for all $s \in [0,S]$.  The estimation of $\cT_{\ep,2}$
is more involved, and the result proved in section \ref{sect:PF1} is
\begin{equation}
  \PP\Big(\sup_{z \in [0,L]} \big|\cT_{\ep,2}(z,s)\big| \ge c
    \ep^{\vartheta} \Big) \le C \ep^{2-\alpha - \vartheta},
  \label{eq:N3}
\end{equation}
with $\vartheta$ satisfying $0 < \vartheta \le p -\alpha < 2
-\alpha$. Then (\ref{eq:C5}), the assumption $\alpha \le 1$ and
Gronwall's inequality give
\begin{equation}
\label{eq:C12}
\PP \Big(\sup_{z \in [0, L]} |q_\ep(z,s)| \ge c\ep^{\vartheta} \Big) \le C
\ep^{2 - \alpha - \vartheta}, \qquad \forall s \in [0,S].
\end{equation}
 
To obtain a uniform estimate in $s$ we note that since the
characteristics cannot intersect, the function $s + q_\ep(z,s)$ is
strictly monotonically increasing in $s$, for all $z \in [0,L]$. More
explicitly,
\begin{equation}
q_\ep(z,s_1) - (s_2-s_1) < q_\ep(z,s_2) < q_\ep(z,s_3) + (s_3-s_2), 
\label{eq:monotone}
\end{equation} 
for all $0 \le s_1 < s_2 < s_3 \le S.$ Let us discretize the
interval $[0,S]$ in steps equal to $\ep^\vartheta S$, with discretization
points $s_j$, for $0 \le j \le S \ep^{-\vartheta}$. For each $s_j$ denote
by $E_j$ the event
\[
E_j = \{ \omega \Big| \sup_{0 \le z \le L} |q_\ep(z,s_j)| \ge c \ep^{\vartheta} \},
\]
with $ \PP(E_j) \le C \ep^{2-\alpha - \vartheta}$ following
from (\ref{eq:C12}). The union of these events has probability
\[
\PP \Big(\cup E_j\Big) \le \sum_{0 \le j \le S\ep^{-\vartheta}}
\PP(E_j) \le C \ep^{2 - \alpha - 2 \vartheta},
\]
where the assumption (\ref{eq:P}) on $p$ guarantees that $ 2 \vartheta
\le 2(p-\alpha) < 2(2 - \alpha).  $ Finally, the monotonicity relation
(\ref{eq:monotone}) gives that
\[
\sup_{z \in [0, L], s \in [0, S]} |q_\ep(z,s)| < c \ep^\vartheta,
\]
on the complement of $\cup E_j$, with probability $1 - O(\ep^{2-\alpha
  - 2 \vartheta})$ and redefined constant $c$. This is the statement
(\ref{lem2.eq2}) in Lemma \ref{thm.01}, with $ \theta = 2 - \alpha - 2
\vartheta$ and $L = c_o T$.

\subsection{Proof of estimate (\ref{eq:N3})}
\label{sect:PF1}
Let us write $\cT_{\ep,2}$ as the sum of two parts
\begin{equation}
  \cT_{\ep,2}(z,s) = s \cJ(z,s) + \ep^{\alpha-2}\cI(z,s),
  \label{eq:NN1}
\end{equation}
where
\begin{align*}
  \cJ(z,s) &= \ep^{1-\alpha}\int_0^z du \, \nu_\ep(u,s),
\\
\cI(z,s) &= \ep^{1-\alpha}\int_0^z du \, \nu_\ep(u,s) \,
\ep^{1-\alpha} \hspace{-0.05in}\int_0^u dv \mu
\Big(\frac{v/c_o}{\ep^\alpha},\frac{v}{\ep^2}\Big).
\end{align*}
We estimate $\cJ$ in section \ref{sect:PF1_1} and $\cI$ in section
\ref{sect:PF1_2}.
\subsubsection{Estimate of $\cJ(z,s)$}
\label{sect:PF1_1}
We begin with proving the basic result
\begin{equation}
\Big| \ep^{1-\alpha} \EE \Big[ \int_{(m-1)\ep}^{m \ep} du \,
  \nu_\ep(u,s) \Big| \mathcal{F}_{\frac{m-1}{\ep}} \Big]\Big| \le
C \ep^{3 -\alpha},  
\label{eq:C11}
\end{equation}
for all $ s \in [0,S]$ and integers $m$ satisfying $1 \le m \le
\lfloor L/\ep\rfloor.$ We
write\footnote{Here we split the integral in three parts to account for all $p$ satisfying (\ref{eq:P}). If
$p$ were $\le 1$, then the part $T_3$ would not be necessary, because 
$(m-1)\ep + \ep^p$ would exceed $m \ep$.}
\begin{equation*}
\EE \Big[\int_{(m-1)\ep}^{(m-1) \ep + \ep^2 } \hspace{-0.1in}du +
  \int_{(m-1)\ep+\ep^2}^{(m-1) \ep + \ep^p}\hspace{-0.1in} du +
  \int_{(m-1)\ep+\ep^p}^{m \ep} \hspace{-0.1in}du \, \nu_\ep(u,s) \Big|
      {\cF_{\frac{m-1}{\ep}}} \Big] = T_1 + T_2 + T_3,
\end{equation*}
with $\nu_\ep$ defined in (\ref{eq:defnue}), for $\zeta_\ep(u,s)$
between $0$ and $s + W_\ep(u-\ep^p) + q_\ep(u-\ep^p,s)$. 
Since
$\nu_\ep$ is bounded, $ |T_1| \le C \ep^2.$ For $T_2$ we can use the
mixing assumption (\ref{eq:asMix}), because $\zeta_\ep(u,s)$ is
$\cF_{\frac{m-1}{\ep}}$ measurable for all $u \le (m-1) \ep + \ep^p$,
and obtain
\begin{align*}
\left|T_2\right| \le C \int_{(m-1)\ep + \ep^2}^{(m-1) \ep + \ep^p}
du \Big(\frac{u}{\ep^2} - \frac{m-1}{\ep}\Big)^{-d} =
\frac{C}{d-1}\ep^2 \left(1-\ep^{(d-1)(2-p)}\right).
\end{align*}
Our assumptions on $p$ and $d$ give that $(d-1)(2-p) > 0$, so $\cT_2$
is $O(\ep^2)$.  For $T_3$ we have $u -\ep^p \ge (m-1) \ep$, so
definition (\ref{eq:defnue}), the tower property of the expectation and
the mixing assumption (\ref{eq:asMix}) give 
\begin{align*}
|T_3| &= c_o^{-1} \Big| \int_{(m-1)\ep + \ep^p}^{m \ep} du \, \EE
\Big[ \EE \Big[ \mu_t \Big( \frac{u/c_o}{\ep^\alpha} +
    \ep^{2-\alpha}\zeta_\ep(u-\ep^p,s), \frac{u}{\ep^2} \Big) \Big|
    \cF_{\frac{u-\ep^p}{\ep^2}} \Big] \Big|
  \cF_{\frac{m-1}{\ep}}\Big] \Big| \\ &\le C \ep
\Big(\frac{\ep^p}{\ep^2}\Big)^{-d} \ll C
\ep^2.
\end{align*}
Here we used that $d > (2-p)^{-1}$ by (\ref{eq:P}), and $d >
2$. Estimate (\ref{eq:C11}) follows.

Now let us define
\[
Y_m = \cJ(m\ep,s) - \cJ((m-1)\ep,s) - \EE \left[\cJ(m\ep,s) -
  \cJ((m-1)\ep,s) \Big| {\cF_{\frac{m-1}{\ep}}} \right],
\]
for\footnote{For convenience we treat henceforth $L/\ep$ as an
  integer. If it is not, there is an additional error term which does
  not make any difference because it tends to zero faster than $\ep$.}
$1 \le m \le L/\ep $, and
\[
X_m = \sum_{j=1}^m Y_j,
\]
which is easily verified to be a martingale, using the tower property
of the conditional expectation. Then,
\begin{equation}
\cJ(m \ep,s) = X_m + \sum_{r=1}^m \EE \Big[\cJ(r\ep,s) -
  \cJ\left( (r-1)\ep,s\right) \Big| \cF_{\frac{r-1}{\ep}}\Big],
\label{eq:E4}
\end{equation}
with terms in the sum bounded like in the left hand side in (\ref{eq:C11}),
and $X_m$ bounded by Doob's  Submartingale Inequality
\[
\PP \Big( \sup_{m \le L/\ep} |X_m| \ge
c\ep^{(2-\alpha+\vartheta)/{2}}\Big) \le
  \frac{\ep^{-(2-\alpha+\vartheta)}}{c^2}\EE \left[ X^2_{L/\ep}\right]
  = \frac{\ep^{-(2-\alpha+\vartheta)}}{c^2}\sum_{j=1}^{L/\ep} \EE
  \left[ Y^2_j\right].
\]
In the last equation we used the tower property of the conditional
expectation,
\begin{align*}
\EE \Big[ \EE \Big[\cJ(r\ep,s) - \cJ((r-1)\ep,s) \Big|
    \cF_{\frac{r-1}{\ep}} \Big] \Big| \cF_{\frac{j}{\ep}} \Big] =
\EE \Big[\cJ(r\ep,s) - \cJ((r-1)\ep,s) \Big|
  \cF_{\frac{j}{\ep}} \Big]
\end{align*}
for $1 \le j < r \le L/\ep$, to get that $ \EE [ Y_j Y_r] = \EE [ Y_j
  \EE [Y_r | \cF_{\frac{j}{\ep}}] ] = 0.$ By the definition of
$\cJ(z,s)$ and the boundedness of $\mu_t$
\begin{equation*}
  \big|\cJ(z,s)-\cJ((m-1)\ep,s)\big| \le C \ep^{2-\alpha}, \qquad
  \forall z \in \big((m-1)\ep,m \ep\big), ~ s \in [0,S],
\end{equation*}
so using (\ref{eq:E4}) terms in the sum estimated by (\ref{eq:C11}) we
get
\begin{equation}
  \sup_{z \in [0,L]}\big|\cJ(z,s)\big| \le C \ep^{2-\alpha} +
  \sup_{m \le L/\ep} \big|X_m\big|.
  \label{eq:E9}
\end{equation}
It remains to bound the martingale $X_m$.

The definition of $Y_m$ and the inequality $(a+b)^2 \le 2 (a^2 + b^2)$
give
\begin{align*}
\EE [ Y_m^2] \le &2 \EE \Big[ |\cJ(m\ep,s) -
  \cJ((m-1)\ep,s)|^2\Big] + \\&2 \EE \Big[ \Big| \EE \Big[
    \cJ(m\ep,s) - \cJ((m-1)\ep,s)\Big|
    \cF_{\frac{m-1}{\ep}}\Big]\Big|^2\Big],
\end{align*}
and using Jensen's inequality
\[
\Big| \EE \Big[ \cJ(m\ep,s) - \cJ((m-1)\ep,s)\Big|
  \cF_{\frac{m-1}{\ep}}\Big]\Big|^2 \le \EE \Big[ \big|\cJ(m\ep,s) -
  \cJ((m-1)\ep,s)\big|^2\Big|\cF_{\frac{m-1}{\ep}}\Big],
\]
and the tower property of the expectation, we get
\begin{align*}
\EE [ Y_m^2] & \le 4 \EE \Big[|\cJ(m\ep,s) - \cJ((m-1)\ep,s)|^2\Big]
\\ & = 4 \ep^{2(1-\alpha)}\EE \Big[\int_{(m-1)\ep}^{m \ep} du
  \,\nu_\ep(u,s)\int_{(m-1)\ep}^{m \ep} dv \, \nu_\ep(v,s)\Big] \\ &=
8 \ep^{2(1-\alpha)} \EE\Big[\int_{(m-1)\ep}^{m \ep} du
  \,\nu_\ep(u,s)\EE \Big[ \int_{u}^{m \ep} dv \,
    \nu_\ep(v,s)\Big|\cF_{\frac{u}{\ep^2}}\Big]\Big].
\label{eq:E10}
\end{align*}
In the last equation we wrote the integral over the square $(u,v) \in
[(m-1)\ep,m\ep] \times [(m-1)\ep,m\ep]$ as twice the integral over its
half, the triangle $u \in [(m-1)]\ep,m\ep]$ and $v \in [u,m\ep]$, and
used the tower property of the expectation. The inner expectation is
estimated as in (\ref{eq:C11}) and the result is
\[
  \EE [ Y_m^2] \le C \ep^{1 + 2(2-\alpha)}.
\]
Substituting in Doob's Submartingale Inequality and then in (\ref{eq:E9}) we get 
\begin{equation}
  \PP\Big( \sup_{z \in [0, L]} |\cJ(z,s)| \ge  c\ep^{(2-\alpha+\vartheta)/2} \Big) \le C
  \ep^{2-\alpha-\vartheta}, \qquad \forall s \in [0,S],
  \label{eq:E13}
\end{equation}
with redefined constants $c$ and $C$.
\subsubsection{Estimate of $\cI(z,s)$}
\label{sect:PF1_2}
We can rewrite the expression of $\cI(z,s)$ as
\begin{equation}
  \cI(z,s) = \ep^{1-\alpha} \int_0^z du\, \nu_\ep(u,s) J(u),
  \label{eq:E14}
\end{equation}
where
\begin{equation}
  J(u) = \ep^{1-\alpha} \int_0^u dv \, \mu
  \Big(\frac{v/c_o}{\ep^\alpha}, \frac{v}{\ep^2}\Big)
  \label{eq:E15}
\end{equation}
can be estimated the same way as $\cJ(z,s)$. We write directly the result
\begin{equation}
  \PP\Big( \sup_{0 \le z \le L} |J(z)| \ge c
  \ep^{(2-\alpha+\vartheta)/2} \Big) \le C \ep^{2-\alpha-\vartheta}.
  \label{eq:E16}
\end{equation}
Moreover, 
\begin{align*}
  \EE \Big[ J^2(z)\Big] &= \ep^{2(1-\alpha)} \iint_0^z du \, d v\, \EE
  \Big[ \mu \Big(\frac{u/c_o}{\ep^\alpha}, \frac{u}{\ep^2}\Big)\mu
    \Big(\frac{v/c_o}{\ep^\alpha}, \frac{v}{\ep^2}\Big) \Big] \\& =
  \ep^{2(1-\alpha)} \iint_0^z du \, d v\, \Phi \Big(\frac{u-v}{c_o
    \ep^\alpha}, \frac{u-v}{\ep^2}\Big),
\end{align*}
and recalling that the autocorrelation $\Phi$ is integrable,
\begin{equation}
  \sup_{0\le z \le L} \EE \Big[ J^2(z)\Big] \le C \ep^{2(2-\alpha)}
  \le C \ep^{2 - \alpha + \vartheta}.
\label{eq:E17}
\end{equation}

Now, we obtain from (\ref{eq:E14}) that
\begin{align*}
  \cI(m\ep,s)-\cI((m-1)\ep,s) = J((m-1)\ep)\,
  \ep^{1-\alpha} \hspace{-0.05in}\int_{(m-1)\ep}^{m
    \ep} \hspace{-0.05in}du \, \nu_\ep(u,s) + \\\ep^{1-\alpha}
  \int_{(m-1)\ep}^{m \ep} du \, \nu_\ep(u,s) \big[
    J(u)-J((m-1)\ep)\big],
\end{align*}
and taking the conditional expectation
\begin{align*}
  \EE \Big[ \cI(m\ep,s)-\cI((m-1)\ep,s) \Big| \cF_{\frac{m-1}{\ep}}
    \Big]=
  J((m-1)\ep)\EE\Big[\ep^{1-\alpha} \hspace{-0.05in}\int_{(m-1)\ep}^{m
      \ep} \hspace{-0.05in}du \, \nu_\ep(u,s) \Big|
    \cF_{\frac{m-1}{\ep}}\Big] + \\\ep^{2(1-\alpha)}
  \EE\Big[\int_{(m-1)\ep}^{m \ep} \hspace{-0.05in}du \, \nu_\ep(u,s)
    \int_{(m-1)\ep}^u \hspace{-0.05in}dv \, \mu
    \Big(\frac{v/c_o}{\ep^\alpha},\frac{v}{\ep^2}\Big)
    \Big|\cF_{\frac{m-1}{\ep}}\Big] .
\end{align*}
We  change the order of integration in the last term 
and condition at $v$ to get
\begin{align*}
  &\EE\Big[\int_{(m-1)\ep}^{m \ep} \hspace{-0.05in}du \, \nu_\ep(u,s)
    \int_{(m-1)\ep}^u \hspace{-0.05in}dv \, \mu
    \Big(\frac{v/c_o}{\ep^\alpha},\frac{v}{\ep^2}\Big)
    \Big|\cF_{\frac{m-1}{\ep}}\Big] = \\& \hspace{0.8in}\EE\Big[\int_{(m-1)\ep}^{m
      \ep} \hspace{-0.05in}dv \, \mu
    \Big(\frac{v/c_o}{\ep^\alpha},\frac{v}{\ep^2}\Big)  \EE \Big[
      \int_v^{m \ep} \hspace{-0.05in} du \, \nu_\ep(u,s) \Big|
      \cF_{\frac{v}{\ep^2}} \Big] \Big|\cF_{\frac{m-1}{\ep}}\Big].
\end{align*}
Then estimates (\ref{eq:C11}) and (\ref{eq:E16}) give 
\begin{align}
  \Big| \EE \Big[ \cI(m\ep,s)-\cI((m-1)\ep,s) \Big|
    \cF_{\frac{m-1}{\ep}} \Big] \Big|1_E &\le C' \ep^{3 -\alpha +
    (2-\alpha+\vartheta)/2} \big[1 + O(\ep^{(2-\alpha-\vartheta)/2})\big]
  \nonumber \\& \le C \ep^{1+(2-\alpha+\vartheta)},
  \label{eq:E18}
\end{align}
where we renamed the constant in the last inequality, used that
$\vartheta < 2 - \alpha$ and $\mu$ is bounded, and let $E$ be the event
with probability $\PP(E) \ge 1 - O(\ep^{2-\alpha-\vartheta})$ so that
$\sup_{z \le L} |J(z)|1_E \le \ep^{(2-\alpha+\vartheta)/2}$. 

We calculate next the second moment of the increments of $\cI(z,s)$,
\begin{align*}
  \EE \Big[ \big[\cI(m\ep,s)-\cI((m-1)\ep,s)\big]^2\Big] \le 2
  \ep^{2(1-\alpha)}\EE \Big[ J^2((m-1)\ep)
    \Big[\int_{(m-1)\ep}^{m\ep} \hspace{-0.05in}du \,
      \nu_\ep(u,s)\Big]^2\Big] + \\ 2 \ep^{4(1-\alpha)}\EE \Big[
    \Big[\int_{(m-1)\ep}^{m\ep} \hspace{-0.05in} du \, \nu_\ep(u,s)
      \int_{(m-1)\ep}^u \hspace{-0.05in} dv \,
      \mu\Big(\frac{v/c_o}{\ep^\alpha},\frac{v}{\ep^2}\Big)\Big]^2\Big]
  = T_1 + T_2.
\end{align*}
To estimate the first term we write the square of the integral of
$\nu_\ep$ as
\begin{align*}
\Big[\int_{(m-1)\ep}^{m\ep} \hspace{-0.05in} du \,
      \nu_\ep(u,s)\Big]^2 = 2   \int_{(m-1)\ep}^{m\ep} \hspace{-0.05in}du \,
\nu_\ep(u,s) \int_u^{m\ep} dv \, \nu_\ep(v,s),
\end{align*} 
condition at $u$, and use (\ref{eq:C11}) and (\ref{eq:E17}) to get 
\begin{align*}
  T_1 &= 4 \ep^{2(1-\alpha)} \EE \Big[ J^2((m-1)\ep)
    \int_{(m-1)\ep}^{m \ep} \hspace{-0.05in} du \, \nu_\ep(u,s) \EE
    \Big[ \int_u^{m \ep} \hspace{-0.05in} dv \, \nu_\ep(v,s) \Big|
      \cF_{\frac{u}{\ep^2}} \Big] \Big] \\ &\le C \ep^{1+ 3(2-\alpha)
    + \vartheta}.
\end{align*}
The second term is
\begin{align*}
  T_2 = 2 \ep^{4(1-\alpha)} \EE \Big[ \int_{(m-1)\ep}^{m
      \ep} \hspace{-0.05in} dv \,
    \mu\Big(\frac{v/c_o}{\ep^\alpha},\frac{v}{\ep^2}\Big)
    \int_{(m-1)\ep}^{m \ep} \hspace{-0.05in} dv' \,
    \mu\Big(\frac{v'/c_o}{\ep^\alpha},\frac{v'}{\ep^2}\Big) \times
    \\ \int_{v}^{m\ep} \hspace{-0.05in} du \, \nu_\ep(u,s)
    \int_{v'}^{m \ep} \hspace{-0.05in} du'\, \nu_\ep(u',s)\Big] =
  T_2^{(i)} + T_2^{(ii)}
\end{align*}
and we split it in two parts
\begin{align*}
  T_2^{(i)} =& 2 \ep^{4(1-\alpha)} \EE \Big[ \int_{(m-1)\ep}^{m
      \ep} \hspace{-0.05in} dv \,
    \mu\Big(\frac{v/c_o}{\ep^\alpha},\frac{v}{\ep^2}\Big)
    \int_{(m-1)\ep}^{v} \hspace{-0.05in} dv' \,
    \mu\Big(\frac{v'/c_o}{\ep^\alpha},\frac{v'}{\ep^2}\Big) \times
    \\ &\Big[ \int_{v}^{m\ep} \hspace{-0.05in} du \, \nu_\ep(u,s)
      \int_{v'}^{v} \hspace{-0.05in} du'\, \nu_\ep(u',s) +
      \int_{v}^{m\ep} \hspace{-0.05in} du \, \nu_\ep(u,s)\int_{v}^{m
        \ep} \hspace{-0.05in} du'\, \nu_\ep(u',s)\big] \Big],
\end{align*}
and
\begin{align*}
  T_2^{(ii)} =& 2 \ep^{4(1-\alpha)} \EE \Big[ \int_{(m-1)\ep}^{m
      \ep} \hspace{-0.05in} dv \,
    \mu\Big(\frac{v/c_o}{\ep^\alpha},\frac{v}{\ep^2}\Big) \int_{v}^{m
      \ep} \hspace{-0.05in} dv' \,
    \mu\Big(\frac{v'/c_o}{\ep^\alpha},\frac{v'}{\ep^2}\Big) \times
    \\&\Big[ \int_{v'}^{m\ep} \hspace{-0.05in} du' \, \nu_\ep(u',s)
      \int_{v}^{v'} \hspace{-0.05in} du\, \nu_\ep(u,s) +
      \int_{v'}^{m\ep} \hspace{-0.05in} du' \, \nu_\ep(u',s)\int_{v'}^{m
        \ep} \hspace{-0.05in} du\, \nu_\ep(u,s)\big] \Big].
\end{align*}
We condition at $v$ in the first term of $T_2^{(i)}$, and use the estimate
(\ref{eq:C11}) for the $u$ integral. The second term involves the
integral
\[
\int_{v}^{m\ep} \hspace{-0.05in} du \, \nu_\ep(u,s)\int_{v}^{m
  \ep} \hspace{-0.05in} du'\, \nu_\ep(u',s) = 2
\int_{v}^{m\ep} \hspace{-0.05in} du \, \nu_\ep(u,s)\int_{u}^{m
  \ep} \hspace{-0.05in} du'\, \nu_\ep(u',s)
 \]
 and we can estimate it using (\ref{eq:C11}), after conditioning at
 $u$. We obtain
 \[
  T_2^{(i)} \le C \ep^{1 + 4(2-\alpha)},
\]
and proceeding in the same manner we obtain a similar bound for $T_2^{(ii)}$.
Gathering the results we have
\begin{equation}
  \EE \Big[ \big[\cI(m\ep,s)-\cI((m-1)\ep,s)\big]^2\Big] \le C \ep^{1
    + 3(2-\alpha) + \vartheta},
  \label{eq:E19}
  \end{equation}
because $\vartheta < (2-\alpha)$.

Now define the martingale $ X_m = \sum_{j=1}^m Y_m,$ with
\[
Y_m = \cI(m\ep,s) - \cI((m-1)\ep,s) - \EE \left[\cI(m\ep,s) -
  \cI((m-1)\ep,s) \Big| {\cF_{\frac{m-1}{\ep}}} \right],
\]
and get the bound
\begin{align}
  |\cI(m \ep,s)|1_E &\le |X_m|1_E + \sum_{j=1}^m \Big|\EE
  \left[\cI(m\ep,s) - \cI((m-1)\ep,s) \Big| {\cF_{\frac{m-1}{\ep}}}
    \right]\Big|1_E \nonumber \\ &\le |X_m|1_E + C m \ep \ep^{2-\alpha
    + \vartheta}
  \label{eq:E20}
\end{align}
where we used (\ref{eq:E18}). The martingale is bounded with Doob's
Submartingale Inequality
\begin{align*}
\PP\Big(\sup_{m \le L/\ep} |X_m| \ge c \ep^{2 -\alpha +
  \vartheta}\Big) \le \frac{\ep^{-2(2-\alpha+\vartheta)}}{c^2} \EE
\Big[ X_{L/\ep}^2\Big] = \frac{\ep^{-2(2-\alpha+\vartheta)}}{c^2}
\sum_{j = 1}^{L/\ep} \EE[Y_j^2] \nonumber \\ \le C \ep^{2-\alpha -
  \vartheta}.
\end{align*}
The equality is implied by the tower property of the expectation, as
in the previous section. The last bound is derived using Jenssen's
inequality and (\ref{eq:E19}).

Finally, definition (\ref{eq:E14}) and (\ref{eq:E16}) give that for any $z \in
\big((m-1)\ep,m\ep\big)$ and $s \in [0,S]$, 
\[
\big|\cI(z,s)-\cI((m-1)\ep,s)\big|1_{E} \le \ep^{2-\alpha} \sup_{u \in
  [0, L]}|J(u)|1_{E} \le \ep^{2-\alpha + (2-\alpha+\vartheta)/2} \le
\ep^{2-\alpha +\vartheta}.
\]
The estimate 
\begin{equation}
\PP\Big(\sup_{0 \le z \le L}|\cI(z,s)| \ge
c\ep^{2-\alpha+\vartheta}\Big) \le C \ep^{2 - \alpha -\vartheta}
\label{eq:E22}
\end{equation}
follows, and together with (\ref{eq:E13}) and definition
(\ref{eq:NN1}) it gives (\ref{eq:N3}).
\subsection{Proof of estimate (\ref{lem2.eq1})}
\label{sect:Pf2}
We introduce the notation
\begin{equation}
\cU_\ep(z,s) = \partial_s \cW_\ep(z/c_o,s),
\label{eq:C14}
\end{equation}
and take the derivative with respect to $s$ in definition
(\ref{eq:Rshift}) to obtain
\begin{equation}
\cU_\ep(z,s) = \ep^{1-\alpha} \int_0^z du \, \gamma_\ep(u,s)
\cU_\ep(u,s) + \ep^{1-\alpha} \int_0^z du \, \gamma_\ep(u,s),
\label{eq:C15}
\end{equation}
with
\begin{equation}
\gamma_\ep(z,s) = c_o^{-1} \mu_t \left( \frac{z/c_o}{\ep^\alpha} +
\ep^{2-\alpha}\left(s + \cW_\ep(z/c_o,s) \right),\frac{z}{\ep^2}\right).
\label{eq:C16}
\end{equation}
The kernel in the first integral satisfies 
\[
\sup_{0 \le u \le L} \ep^{1-\alpha} |\gamma_\ep(u,s)| \le C
\ep^{1-\alpha}, 
\]
  because $\mu_t$ is bounded. The second integral is
estimated exactly like in section \ref{sect:PF1_1}, once we we prove
the equivalent of the basic estimate (\ref{eq:C11}) for $\gamma_\ep$.

Let us rewrite (\ref{eq:C16}) as
\begin{align*}
\gamma_\ep(u,s) = c_o^{-1} \Big[\mu_t \Big( \frac{u/c_o}{\ep^\alpha} +
  \ep^{2-\alpha}\left(s+\cW_\ep( (u-\ep^p)/{c_o},s ) \right),\frac{u}{\ep^2}\Big)
  + \cE_\ep(u,s)\Big],
\end{align*}
with $p$ as in (\ref{eq:P}) and error
\begin{align*}
\cE_\ep(u,s) =& \mu_t \left( \frac{u/c_o}{\ep^\alpha} +
\ep^{2-\alpha}\left(s+\cW_\ep(u/c_o,s)\right),\frac{u}{\ep^2}\right) -
\nonumber \\ &\mu_t \left( \frac{u/c_o}{\ep^\alpha} +
\ep^{2-\alpha}\left(s+\cW_\ep
((u-\ep^p)/c_o,s)\right),\frac{u}{\ep^2}\right).
\end{align*}
We have
\begin{equation}
\ep^{1-\alpha} \EE\Big[ \int_{(m-1)\ep}^{m \ep} du \, \gamma_\ep(u,s)
  \Big| \cF_{\frac{m-1}{\ep}} \Big] = c_o^{-1} \ep^{1-\alpha}
\left(T_1 + T_2 + T_3\right),
\label{eq:C19}
\end{equation}
where
\begin{align*}
T_1 &= \EE\Big[ \int_{(m-1)\ep}^{m \ep} \hspace{-0.05in}du \,
  \cE_\ep(u,s) \Big| \cF_{\frac{m-1}{\ep}} \Big], \\ T_2 &= \EE\Big[
    \int_{(m-1)\ep}^{(m-1) \ep + \ep^2 } \hspace{-0.1in}du \, \mu_t
    \Big( \frac{u/c_o}{\ep^\alpha} + \ep^{2-\alpha}\big(s +
    \cW_\ep((u-\ep^p)/c_o,s) \big) ,\frac{u}{\ep^2}\Big) \Big|
    \cF_{\frac{m-1}{\ep}} \Big], \\ T_3 &= \EE\Big[
    \int_{(m-1)\ep+\ep^2}^{m\ep} \hspace{-0.1in}du \,\mu_t \Big(
    \frac{u/c_o}{\ep^\alpha} + \ep^{2-\alpha}\big(s +
    \cW_\ep((u-\ep^p)/c_o,s)\big) ,\frac{u}{\ep^2}\Big) \Big|
    \cF_{\frac{m-1}{\ep}} \Big].
\end{align*}
Definition (\ref{eq:C14}), the mean value theorem, the boundedness of
$\mu_{tt}$ and the estimate
\begin{align*}
\left|\cW_\ep(u/c_o,s)-\cW_\ep((u-\ep^p)/c_o,s)\right|\le C\ep^{-1+p} ,
\end{align*}
give that $ \left| T_1 \right| \le C \ep^{2 + p-\alpha} \le C \ep^2.
$ The last inequality is because $p > \alpha.$ The estimate $ \left|
T_2 \right| \le C \ep^2 $ follows from the boundedness of $\mu_t$, and
$T_3 = T_{3}^{(i)} + T_{3}^{(ii)},$ where
\begin{align*}
T_{3}^{(i)} &= \EE \Big[ \int_{(m-1)\ep+\ep^2}^{(m-1)\ep + \ep^p} du \,
  \mu_t \Big(\frac{u/c_o}{\ep^{\alpha}} + \ep^{2-\alpha}(s +
  \cW_\ep((u-\ep^p)/c_o,s)),\frac{u}{\ep^2}\Big)\Big|
  \cF_{\frac{m-1}{\ep}}\Big],\\
T_{3}^{(ii)} &= \EE \Big[ \int_{(m-1)\ep+\ep^p}^{m \ep} du \, \mu_t
  \Big(\frac{u/c_o}{\ep^{\alpha}} + \ep^{2-\alpha}\big(s +
  \cW_\ep((u-\ep^p)/c_o,s)\big) ,\frac{u}{\ep^2}\Big)\Big|
  \cF_{\frac{m-1}{\ep}}\Big].
\end{align*}
In $T_{3}^{(i)}$ we note that $\cW_\ep((u-\ep^p)/c_0,s)$ is $\cF_{\frac{m-1}{\ep}}$
measurable, and use the mixing assumption (\ref{eq:asMix})
\[
\left|T_{3}^{(i)}\right| \le C \hspace{-0.05in}\int_{(m-1)\ep +
  \ep^2}^{(m -1)\ep + \ep^p} \hspace{-0.05in}d u \left(
\frac{u}{\ep^2} - \frac{m-1}{\ep}\right)^{-d} \hspace{-0.1in}= \frac{C
  \ep^2}{d-1} \left[1 - \ep^{(d-1)(2-p)}\right] = O(\ep^2).
\]
In $T_{3}^{(ii)}$ we have $u - \ep^p \ge (m-1)\ep$, so we use the
tower property of the expectation and the mixing assumption
(\ref{eq:asMix}) to obtain
\begin{align*}
|T_{3}^{(ii)}| = &\Big|\hspace{-0.02in}\int_{(m-1)\ep + \ep^p}^{m
  \ep} \hspace{-0.25in}du \, \EE \Big[ \EE \Big[
    \mu_t\Big(\frac{u/c_o}{\ep^\alpha} + \ep^{2-\alpha}\big(s +
    \cW_\ep((u-\ep^p)/c_o,s)\big), \frac{u}{\ep^2}\Big)\Big|
    \cF_{\frac{u-\ep^p}{\ep^2}} \Big] \Big| \cF_{\frac{m-1}{\ep}}\Big]
\Big| \nonumber \\&\le C \ep \Big(\frac{\ep^p}{\ep^2}\Big)^{-d} \ll C
\ep^2.
\end{align*}
The last inequality is because $d > (2-p)^{-1}$.  Gathering the
results we get the equivalent of (\ref{eq:C11})
\[
\Big| \ep^{1-\alpha} \EE \Big[ \int_{(m-1)\ep}^{m \ep} du \,
  \gamma_\ep(u,s) \Big| \mathcal{F}_{\frac{m-1}{\ep}} \Big]\Big| \le
C \ep^{3 -\alpha},
\]
and derive as in section \ref{sect:PF1_1} that 
\[
\PP\Big(\sup_{z \in [0, L]} \Big|\ep^{1-\alpha} \int_0^z du\,
\gamma_\ep(u,s)\Big| \ge c \ep^{(2-\alpha+\vartheta)/2} \Big) \le C \ep^{2-\alpha-\vartheta}.
\]
By assumption $\vartheta < 2-\alpha$, so we conclude from
(\ref{eq:C15}), the assumption $\alpha \le 1$ and Gronwall's
inequality that
\begin{equation}
\PP\Big(\sup_{z \in [0, L]} \big|\cU_\ep(z,s)\big| \ge c
\ep^\vartheta\Big) \le C \ep^{2-\alpha - \vartheta}, \qquad \forall s
\in [0,S].
\label{eq:E32}
\end{equation}

To make the result uniform in $s$, we now estimate the variation of
$\cU_\ep$, denoted by
\[
\Delta \cU_\ep(z,s,\Delta s) = \cU_\ep(z,s+\Delta s)
-\cU_\ep(z,s), 
\]
for all $\Delta s$ satisfying $0 \le \Delta s \le \ep^{\vartheta} S$. 
Equation (\ref{eq:C15}) gives 
\begin{align}
\Delta \cU_\ep(z,s,\Delta s) =
\ep^{1-\alpha} \hspace{-0.05in}\int_0^z \hspace{-0.05in} du \, \Delta
\gamma_\ep(u,s ,\Delta s) \left[ 1 + \cU_\ep(u,s)\right] + \nonumber
\\ \ep^{1-\alpha} \int_0^z du \, \gamma_\ep(u,s+\Delta s) \Delta
\cU_\ep(u,s,\Delta s),\qquad \label{eq:C24p}
\end{align}
where we introduced the notation $ \Delta \gamma_\ep(z,s,\Delta s) =
\gamma_\ep(z,s+\Delta s) - \gamma_\ep(z,s).  $ The kernel multiplying
$\Delta \cU_\ep$ in the last integral in (\ref{eq:C24p}) is uniformly
bounded by
\[
\sup_{u \in [0,L] , \Delta s \in [0, \ep^\vartheta S]}
\Big|\ep^{1-\alpha} \gamma_\ep(u,s+\Delta s)\Big| \le C
\ep^{1-\alpha},
\]
  because $\mu_t$ and therefore $\gamma_\ep$ are bounded.
Moreover, using definition (\ref{eq:C16}) and the mean value theorem
in the definition of $\Delta \gamma_\ep$ we get
\begin{align*}
\Delta \gamma_\ep(z,s,\Delta s) = c_o^{-1} \ep^{2-\alpha} \mu_{tt}
\Big(\frac{z/c_o}{\ep^{\alpha}}+\ep^{2-\alpha}
\zeta_\ep(z,s,\Delta s),\frac{z}{\ep^2}\Big)\times \nonumber \\\Big[
  \Delta s + \cW_\ep(z/c_o,s+\Delta s) - \cW_\ep(z/c_o,s)\Big],
\end{align*}
with $\zeta_\ep(z,s,\Delta s)$ somewhere between $s +
\cW_\ep(z/c_o,s)$ and $s + \Delta s + \cW_\ep(z/c_o,s + \Delta s)$. By
assumption $\mu_{tt}$ is uniformly bounded. To bound
the second factor, we recall the monotonicity relation
(\ref{eq:monotone}), and obtain that
\begin{equation*}
  0 \le \Delta s + \cW_\ep(z/c_o,s+\Delta s) - \cW_\ep(z/c_o,s) \le
  \ep^\vartheta S + \cW_\ep(z/c_o,s+\ep^\vartheta S) -
  \cW_\ep(z/c_o,s), 
\end{equation*}
with right hand side estimated in (\ref{eq:C12}) for any given $s$ and
$s + \ep^\vartheta S$. Consequently,
\begin{equation*}
\PP \Big( \sup_{z \in [0,L], \Delta s \in [0, \ep^\vartheta
    S]}\big|\Delta \gamma_\ep(z,s,\Delta s)\big| \ge c\ep^{ 2-\alpha +
  \vartheta} \Big) \le C \ep^{2 - \alpha - \vartheta},
\end{equation*}
and recalling the bound (\ref{eq:E32}), we get 
\begin{align*}
\PP \Big( \sup_{z \in [0, L], \Delta s \in [0, \ep^\vartheta
    S]}\Big|\ep^{1-\alpha} \hspace{-0.05in} \int_0^z du \Delta
\gamma_\ep(u,s,\Delta s)\big[1+\cU_\ep(u,s)\big]\Big| \ge c \ep^{3
  -2\alpha + \vartheta} \Big) \le C \ep^{2 - \alpha - \vartheta}.
\end{align*}
The estimate of $\Delta \cU_\ep$ follows by Gronwall's inequality from
this result and equation (\ref{eq:C24p})
\begin{equation}
\PP \Big( \sup_{0\le z \le L, 0 \le \Delta s \le \ep^\vartheta S}\Big|
\Delta \cU_\ep(z,s,\Delta s)\Big| \ge c \ep^{3-2\alpha+\vartheta} \Big) \le C
\ep^{2-\alpha - \vartheta}.
\label{eq:E40}
\end{equation}

Now we can bound $\cU_\ep(z,s)$ for all $s \in (s_{j-1},s_j)$, with
$s_j = j \ep^\vartheta S$,  using the triangle inequality
\[
\big|\cU_\ep(z,s)\big| \le \big|\cU_\ep(z,s_{j-1})\big|+\big|\Delta
\cU_\ep(z,s_{j-1},s-s_{j-1})\big|.
\]
We obtain from (\ref{eq:E32}), (\ref{eq:E40}) and the assumption
$\alpha \le 1$ that
\begin{equation*}
\sup_{z \in [0,L], s \in [s_{j-1}, s_{j}]}\big|\cU_\ep(z,s)\big|
1_{E_j} < c \ep^\vartheta.
\end{equation*}
Here $E_j$ is an event with probability $\PP(E_j) \ge 1 -
O(\ep^{2-\alpha - \vartheta})$. There are $\ep^{-\vartheta}$ such events and
on their intersection $E = \cap_{j} E_j$, satisfying
$\PP(E) \ge 1 - O(\ep^{2-\alpha - 2 \vartheta})$,
\begin{equation*}
\sup_{z \in [0,L], s \in [0,S]}\big|\cU_\ep(z,s)\big| 1_{E} < c
\ep^\vartheta.
\end{equation*}
This is the statement (\ref{lem2.eq1}). 

\subsection{Regime 2}
\label{sect:PfReg2}
The case $\beta \le 1$, $\alpha = 2$, where 
\begin{equation}
\label{eq:B10}
\cW_\ep(\tau,s) = \frac{1}{\ep} \int_0^\tau d u \, \mu
\left(\frac{u}{\ep^2} + s + \cW_\ep(u,s),\frac{c_o
  u}{\ep^\beta}\right),
\end{equation}
can be handled with the change of variables
\begin{equation}
\zeta = u + \ep^{2}s + \ep^2 \cW_\ep(u,s),
\label{eq:B11}
\end{equation}
so that
\begin{equation}
d \zeta = \left[1 + \ep \mu \left(\frac{u}{\ep^2} + s +
  \cW_\ep(u,s),\frac{c_o u}{\ep^\beta}\right)\right] du.
\end{equation}
Since $\mu$ is bounded by assumption, we know from the implicit
function theorem that we can invert (\ref{eq:B11}) to get a unique
function $u = u(\zeta)$, and with, $u = \zeta + O(\ep^2)$.
Substituting (\ref{eq:B11}) in (\ref{eq:B10}), we obtain that
\begin{align*}
\cW_\ep(\tau,s) &\approx \frac{1}{\ep} \int_0^\tau d\zeta \,
\frac{\mu\left(\frac{\zeta}{\ep^2} , \frac{c_o \zeta}{\ep^\beta} -
  \ep^{2-\beta} c_o(s+\cW_\ep(\zeta,s))\right)}{1 + \ep
  \mu\left(\frac{\zeta}{\ep^2} , \frac{c_o \zeta}{\ep^\beta} -
  \ep^{2-\beta} c_o(s+\cW_\ep(\zeta,s))\right)},
\end{align*}
and expanding the right hand side in $\ep$,
\begin{align*}
\cW_\ep(\tau,s) &\approx
\frac{1}{\ep} \int_0^\tau d\zeta \, \mu\left(\frac{\zeta}{\ep^2} ,
\frac{c_o \zeta}{\ep^\beta} - \ep^{2-\beta} c_o(s+\cW_\ep(\zeta,s))\right)
-\int_0^\tau d\zeta \,
\mu^2\left(\frac{\zeta}{\ep^2} , \frac{c_o \zeta}{\ep^\beta}\right).
\end{align*}
Here we used that $\mu$ and $\mu_z$ are bounded, and the approximate sign 
stands for equal up to absolute errors that
tend to zero as $\ep \to 0$.  Now we can use the same analysis as in sections \ref{sect:Pf1}-\ref{sect:Pf2} 
to show that we can write $\cW_\ep(\tau,s)$  as stated in Lemma \ref{thm.01}.

\section{Proof of Lemma \ref{lem.apD}}
\label{ap:D}
We wish to show that
\begin{align}
\cT^\ep(\tau;s,y) = \int_0^{\tau} d u \Big[ \mu_z
  \left(\frac{u}{\ep^\alpha} + \ep^{2-\alpha}\Big(\frac{s+y}{2} +
  \cW_\ep(u,s)\Big), \frac{c_o u}{\ep^2} +  \frac{c_o(s-y)}{2} \right)
  \times \nonumber \\ \mu_z \left(\frac{u}{\ep^\alpha}+
  \ep^{2-\alpha}\big(s + \cW_\ep(u,s)\big), \frac{c_o u}{\ep^2}
  \right)+ \partial_z^2 \Phi\left(0,\frac{c_o(s-y)}{2} \right)
  \Big] \label{eq:D0}
\end{align}
converges to $0$ as $\ep \to 0$, uniformly in $\tau
\le T$ and $0 \leq y \leq  s \in [0,S]$.  
We write (\ref{eq:D0}) as the sum $\cT^\ep
= \cT_1^\ep + \cT_2^\ep$, where the first term is given by
\begin{align*}
& \cT_1^\ep(\tau) = \int_0^{\tau} d u \Big[ \mu_z
    \left(\frac{u}{\ep^\alpha} + \ep^{2-\alpha}\Big(\frac{s+y}{2} +
    \cW_\ep(u,s)\Big), \frac{c_o u}{\ep^2} + \frac{c_o(s-y)}{2}
    \right) \times \\ &\hspace{0.2in}\mu_z \left(\frac{u}{\ep^\alpha}+
    \ep^{2-\alpha}\big(s + \cW_\ep(u,s)\big), \frac{c_o u}{\ep^2}
    \right) -\mu_z \left(\frac{u}{\ep^\alpha}+ \ep^{2-\alpha}\big(s +
    \cW_\ep(u-\ep^p,s)\big), \frac{c_o u}{\ep^2} \right) \times
    \\ &\hspace{0.2in}\mu_z \left(\frac{u}{\ep^\alpha} +
    \ep^{2-\alpha}\Big(\frac{s+y}{2} + \cW_\ep(u-\ep^p,s)\Big),
    \frac{c_o u}{\ep^2}+ \frac{c_o(s-y)}{2} \right) \Big],
\end{align*}
and the second term is 
\begin{align*}
\cT_2^\ep(\tau) = \int_0^{\tau} d u \Big[ \mu_z
\left(\frac{u}{\ep^\alpha} + \ep^{2-\alpha}\Big(\frac{s+y}{2} +
\cW_\ep(u-\ep^p,s)\Big), \frac{c_o u}{\ep^2} + \frac{c_o(s-y)}{2}
\right) \times \\ \mu_z \left(\frac{u}{\ep^\alpha}+
\ep^{2-\alpha}\big(s + \cW_\ep(u-\ep^p,s)\big), \frac{c_o
  u}{\ep^\beta} \right)+ \partial_z^2 \Phi\left(0,\frac{c_o(s-y)}{2}
\right)\Big],
\end{align*}
with $p$ as in (\ref{eq:P}). We suppressed $s$ and $y$ in the
arguments to simplify notation. The mean value theorem and the
definition of $\cW_\ep$ allows us to bound $\cT_1^\ep$ as
\begin{equation} 
|\cT_1^\ep| \le C \ep^{p+1-\alpha},
\label{eq:D1}
\end{equation}
with constant $C$ that is independent of $\ep$, $\tau,
s $ and $y$. 

To estimate $\cT_2^\ep$ we use an argument that is similar to that
used in Appendix \ref{sect:PF1}.  We begin with
\begin{align}
\label{eq:D2}
\EE \left[ \cT_2^\ep(m \ep) - \cT_2^\ep((m-1)\ep)\Big|
  \cF_{\frac{m-1}{\ep}} \right] = \sum_{j=1}^3 J_j,
\end{align}
which we write as the sum of three terms.  The first satisfies
\begin{equation*}
\left|J_1\right| = \left|\EE \left[ \cT_2^\ep((m-1) \ep + \ep^2) -
  \cT_2^\ep((m-1)\ep)\Big| \cF_{\frac{m-1}{\ep}} \right]\right| \le C
\ep^2,
\end{equation*}
because $\mu_z$ is bounded. The second term is
\begin{align*}
J_2 = \EE \left[ \cT_2^\ep((m-1) \ep+\ep^p) -
  \cT_2^\ep((m-1)\ep+\ep^2)\Big| \cF_{\frac{m-1}{\ep}} \right],
\end{align*}
and we can bound it using the mixing assumption (\ref{eq:asMix2}).
Since $\cW_\ep(u-\ep^p,s)$ is $\cF_{\frac{m-1}{\ep}}$ measurable for
$u \le (m-1)\ep + \ep^p$, we obtain from (\ref{eq:asMix2}) and the
definition of $\cT_2^\ep$ that
\begin{align*}
\left|J_2 - \int_{(m-1)\ep + \ep^2}^{(m-1)\ep + \ep^p}
\hspace{-0.05in}d y \left[ \partial_z^2 \Phi
  \left(0,\frac{c_o(s-y)}{2} \right) - \partial_z^2 \Phi
  \left(\frac{\ep^{2-\alpha}(y-s)}{2},\frac{c_o(s-y)}{2} \right)
  \right] \right| \le \nonumber \\ C \int_{(m-1)\ep + \ep^2}^{(m-1)\ep
  + \ep^p} \left( \frac{u}{\ep^2} - \frac{m-1}{\ep} \right)^{-d} =
\frac{C \ep^2}{d-1}\left(1 - \ep^{(d-1)(2-p)}\right).
\end{align*}
The integrand in the left hand side is $O(\ep^{2-\alpha})$, and since
$p + 2-\alpha > 2$ we get $ \left|J_2 \right| \le C \ep^{2}.  $ In
the last term
\begin{align*}
J_3 = \EE \left[ \cT_2^\ep(m\ep) - \cT_2^\ep((m-1)\ep+\ep^p)\Big|
  \cF_{\frac{m-1}{\ep}} \right]
\end{align*}
we have $u - \ep^p \ge (m-1) \ep$ and use the tower property of
the expectation to write 
\begin{align*}
J_3 = \EE \left[ \int_{(m-1)\ep + \ep^p}^{m \ep} \hspace{-0.1in}d u
  \left\{ \EE \left[ \mu_z \Big(\frac{u}{\ep^\alpha} +
    \ep^{2-\alpha}\big(\frac{s+y}{2} + \cW_\ep(u-\ep^p,s)\big),
    \frac{c_o u}{\ep^2} + \frac{c_o(s-y)}{2} \Big)
    \right. \right. \right. \hspace{-0.06in}\times
    \\ \left. \left. \left. \mu_z \left( \frac{u}{\ep^\alpha}+
    \ep^{2-\alpha}\big(s + \cW_\ep(u-\ep^p,s)\big), \frac{c_o
      u}{\ep^2} \right) \Big| \cF_{\frac{u-\ep^p}{\ep^2}}\right]+
  \partial_z^2 \Phi \left(0,\frac{c_o(s-y)}{2} \right) \right\} \Big|
  \cF_{\frac{m-1}{\ep}} \right],
\end{align*}
The inner expectation equals  
 \[ -\partial_z^2
\Phi\left(\frac{\ep^{2-\alpha}(y-s)}{2},\frac{c_o (s-y)}{2}\right) +
O(\ep^{(2-p)d}) = -\partial_z^2\Phi\left(0,\frac{c_o (y-s)}{2}\right)
+ O(\ep^{2-\alpha}) ,\] because $ (2-p)d > 2 -\alpha.  $
Thus, $ |J_3| \le C \ep^2, $ and from (\ref{eq:D2}) we find
\begin{equation}
\label{eq:D4}
\left| \EE \left[ \cT_2^\ep(m \ep) - \cT_2^\ep((m-1)\ep)\Big|
  \cF_{\frac{m-1}{\ep}} \right] \right| \le C \ep^2. 
\end{equation}

Now let
\begin{equation}
Y_m = \cT_2^\ep(m \ep) - \cT_2^\ep((m-1)\ep) - \EE \left[ \cT_2^\ep(m
  \ep) - \cT_2^\ep((m-1)\ep)\Big| \cF_{\frac{m-1}{\ep}} \right]
\end{equation}
and define the martingale $ X_m = \sum_{j=1}^m Y_j.  $ We have
\begin{align*}
\Big| \cT_2^\ep(m \ep) \Big| &\le |X_m| + \sum_{j=1}^m \Big|\EE
\Big[ \cT_2^\ep(m \ep) - \cT_2^\ep((m-1)\ep)\Big|
  \cF_{\frac{m-1}{\ep}} \Big] \Big|\le |X_m| + C m
\ep^2,
\end{align*}
with the second bound following from (\ref{eq:D4}). We also get from
the definition of $\cT_2^\ep$ and the boundedness of $\mu_z$ that
\[
\Big| \cT_2^\ep(\tau)-\cT_2^\ep(m \ep) \Big| \le C \ep, \qquad 
\forall \, \tau \in (m \ep, (m+1) \ep),
\]
and therefore
\begin{equation}
\label{eq:D8}
\sup_{\tau \in [0, T]} |\cT_2^\ep(\tau)| \le \sup_{m \le T/ \ep} |X_m|
+ C \ep.
\end{equation}

It remains to bound $X_m$, using the Submartingale Inequality
\begin{equation}
\label{eq:D9}
\PP \left(\sup_{m \le T /\ep } |X_m| \ge c \ep^{\frac{2-\theta}{2}} \right) \le
\frac{\ep^{\theta-2}}{c^2} \EE \left[X_{T/\ep}^2\right],
\end{equation}
for $\theta$ as in Lemma \ref{thm.01}, 
where we suppose for simplicity that $T/\ep$ is integer.
 We have
\begin{equation}
\label{eq:D9B}
\EE\left[X_{T /\ep}^2\right] = \sum_{j=1}^{T/\ep}\EE\left[Y_j^2\right]
\le 4 \sum_{j=1}^{T /\ep } \EE \left[ \left|
  \cT_2^\ep(j\ep)-\cT_2^\ep((j-1)\ep)\right|^2\right],
\end{equation}
where we used the tower property of the expectation and Jenssen's
inequality. To calculate the last expectation we denote in short by
$n_\ep(u)$ the integrand in the definition of $\cT_2^\ep$. Again, we
suppress $s$ and $y$ for simplicity. We have
\begin{align*}
\EE \left[ \iint_{(m-1)\ep}^{m \ep} d u \, dv \, n_\ep(u) n_\ep(v)
  \right] &= 2 \EE \left[ \int_{(m-1)\ep}^{m \ep} d u \, n_\ep(u)
  \int_{u}^{m \ep } dv \, n_\ep(v) \right] \\ & = 2 \EE \left[
  \int_{(m-1)\ep}^{m \ep} d u \, n_\ep(u) \EE \left[\int_{u}^{m \ep }
    dv \, n_\ep(v) \Big| \cF_{\frac{u}{\ep^2}} \right] \right], 
\end{align*}
with the conditional expectation bounded as in (\ref{eq:D4}). Thus, 
\[
\EE\left[ \left| \cT_2^\ep(j\ep)-\cT_2^\ep((j-1)\ep)\right|^2\right]
\le C \ep^3,
\]
and substituting in (\ref{eq:D9})-(\ref{eq:D9B}) we get the
probabilistic bound on $X_m$. The estimate (\ref{eq:D8}) becomes
\begin{equation}
\label{eq:D12}
\PP \left(\sup_{\tau \in [0, T]} |\cT_2^\ep(\tau)| \ge c
\ep^{\frac{2-\theta}{2}}\right) \le CT \ep^{\theta}.
\end{equation}

Results (\ref{eq:D1}) and (\ref{eq:D12}) give that $\cT^\ep(\tau;s,y)$
converges in probability to $0$, uniformly in $\tau$ and pointwise in
$s$ and $y$. More explicitly, assumption (\ref{eq:P}) gives  
$p + 1- \alpha > 1 > (2-\theta)/2$, so 
$\cT_2^\ep$ dominates $\cT_1^\ep$ and 
\begin{equation}
\label{eq:cor}
\PP \left(\sup_{\tau \in [0, T]} |\cT^\ep(\tau;s,y)| \ge c
\ep^{\frac{2-\theta}{2}}\right) \le CT \ep^{\theta},
\end{equation}
for all $s$ and $y$ satisfying $0 \le y \le s \le S$.  
Next  we obtain from definition
(\ref{eq:D0}) and the assumptions on $\mu$ and $\Phi$ that  
\begin{equation}
\sup_{\tau \in [0,T], 0 \le y \le s \le S} |\partial_y \cT^\ep(\tau;s,y)|
\le C .
\label{eq:cor1}
\end{equation}
The derivative $\partial_s \cT^\ep$
involves $\partial_s \cW_\ep$, for which we have the result in Lemma
\ref{thm.01}, so 
\begin{equation}
\PP \Big( \sup_{\tau \in [0,T], 0 \le y \le s \le S} |\partial_s
\cT^\ep(\tau;s,y)| \ge c\Big) \le C \ep^{\theta}. 
\label{eq:cor2}
\end{equation}

 We  discretize the interval $[0, S]$ in uniform
steps $\ep^{r}$ with $r$ to be specified below,    and obtain
\begin{equation}
\PP \Big( \sup_{j\leq m \leq S/\ep^{r}, \tau \in [0,T] } | 
\cT^\ep(\tau;m \ep^r, j\ep^r)| \ge c\ep^{\frac{2-\theta}{2}} \Big) \le CT \ep^{\theta-2r}. 
\label{eq:cor2b}
\end{equation}
The factor $\ep^{-2r}$ in the   bound is because we have $\ep^{2 r}$ terms 
enumerated by $(m,j)$. 
In view of (\ref{eq:cor1}) and  (\ref{eq:cor2}) we then find
\begin{equation}
\PP \Big( \sup_{y\leq s \leq S, \tau \in [0,T] } | 
\cT^\ep(\tau; s, y)| \ge 
c \big(\ep^{\frac{2-\theta}{2}}+\ep^r \big)    \Big) \le C T\ep^{\theta-2r}.
\label{eq:cor2c}
\end{equation}
Moreover, letting $T \mapsto T\ep^q$, and absorbing the new $T$ in the constant $C$, we get 
 \begin{equation}
\PP \Big( \sup_{y\leq s \leq S, \tau \in [0,T\ep^q] } | 
\cT^\ep(\tau; s, y)| \ge 
c  \big(\ep^{\frac{2-\theta}{2} }+\ep^r\big) \Big) \le C \ep^{\theta-2r+q}.
\label{eq:cor2f}
\end{equation}
Finally, let $r = (\theta-\vartheta)/2$, for some $\vartheta$ satisfying 
$\max[0,2(\theta-1)] < \vartheta <  \theta$, so that  
\[
r = \frac{\theta-\vartheta}{2} < \frac{2-\theta}{2},
\]
and therefore
 \begin{equation}
\PP \Big( \sup_{y\leq s \leq S, \tau \in [0,T\ep^q] } | 
\cT^\ep(\tau; s, y)| \ge 
c \ep^{\frac{\theta-\vartheta}{2}} \Big) \le C \ep^{\vartheta +q}.
\label{eq:cor2g}
\end{equation}
For $\varphi$ a smooth function this then gives 
  Lemma \ref{lem.apD}.  NB any starting point

\bibliographystyle{siam} \bibliography{LTDep}

\end{document}